 \documentclass[12pt]{article}

\usepackage{amsmath}
\usepackage{amsfonts}
\usepackage{amssymb}
\input xy
\xyoption{all}
\usepackage{srcltx}

\usepackage{hyperref}

\usepackage{xcolor}

\usepackage[normalem]{ulem}

\usepackage{marginnote}

\input xy
\xyoption{all}

\newtheorem{theo}{\bf \thesection.\arabic{abz}. Theorem }

\newtheorem{lemm}{\bf \thesection.\arabic{abz}. Lemma}

\newtheorem{coro}{\bf \thesection.\arabic{abz}. Corollary}

\newtheorem{defi}{\bf  \thesection.\arabic{abz}. Definition }

\newtheorem{exa}{\bf \thesection.\arabic{abz}. Example}

\newtheorem{rema}{\bf \thesection.\arabic{abz}. Remark}

\newcounter{abz}[section]
\newcounter{equ}[section]
\newcounter{equu}[section]
\newcommand{\abz}{\refstepcounter{abz}}
\newcommand{\equ}{\refstepcounter{equ}}

\def\ad{\mathrm{ad\,}}

\def\codim{\mathrm{codim}}

\def\diag{\mathrm{diag}}

\def\Id{\mathrm{Id}}

\def\F{\mathcal{F}}

\def\G{\Gamma}

\def\Pf{\mathrm{Pf}}

\def\rk{\mathop{\mathrm{rank}}}

\def\Tr{\mathop{\mathrm{Tr}}}

\def\V{\mathcal{V}}

\def\K{\mathbb{K}}

\def\R{\mathbb{R}}

\def\P{\mathbb{P}} \def\C{\mathbb{C}}

\def\T{\Theta}

\def\al{\alpha}
\def\be{\beta}

\def\ga{\gamma}

\def\d{\partial}

\def\la{\lambda}

\def\om{\omega}

\newcommand{\iprod}{\mathbin{\lrcorner}}
\def\qed{$\square$}

\title{Webs, Nijenhuis operators, and heavenly PDEs}
\author{
\and
Andriy Panasyuk\\
Faculty of Mathematics and Natural Sciences\\
 Cardinal Wyszy\'{n}ski
University\\
Wóycickiego 1/3, 01-938 Warsaw, Poland\\
\\
Adam Szereszewski\\
Faculty of Physics, University of Warsaw\\
Pasteura 5,  02-093 Warsaw, Poland\\
Adam.Szereszewski@fuw.edu.pl
}
\date{}

\providecommand{\bysame}{\leavevmode\hbox to3em{\hrulefill}\thinspace}
\providecommand{\MR}{\relax\ifhmode\unskip\space\fi MR }
\providecommand{\MRhref}[2]{%
  \href{http://www.ams.org/mathscinet-getitem?mr=#1}{#2}
}
\providecommand{\href}[2]{#2}

\begin{document}
\bibliographystyle{amsalpha}

\maketitle

\begin{abstract}
 In 1989 Mason and Newman proved that there is a 1-1-correspondence between self-dual metrics satisfying Einstein vacuum equation (in complex case or in neutral signature) and pairs of commuting parameter depending vector fields $X_1(\lambda),X_2(\lambda)$ which are divergence free with respect to some volume form. Earlier (in 1975) Plebański showed instances of such vector fields depending of one function of four variables satisfying the so-called I or II Plebański heavenly PDEs.  Other PDEs leading to Mason--Newman vector fields are also known in the literature: Husain--Park (1992--94), {\color{black}Grant (1993), }Schief (1996). In this paper we discuss these matters in the context of the web theory, i.e. theory of collections of foliations on a manifold, understood from the point of view of Nijenhuis operators. In particular we  show how to apply this theory for constructing new ``heavenly'' PDEs based on different normal forms of Nijenhuis operators in 4D, which are integrable similarly to their predecessors. Relation with the Hirota dispersionless systems of PDEs and the corresponding Veronese webs, which was recently observed by  Konopelchenko--Schief--Szereszewski,  is established in all the cases. We also discuss some higher dimensional generalizations of the ``heavelny'' PDEs and the existence of related vacuum Einstein metrics  in 4D-case.
\end{abstract}

\tableofcontents
\section*{Introduction}

In 1989 Mason and Newman \cite{masonNewman} proved that there is a 1-1-correspondence between self-dual metrics satisfying Einstein vacuum equation (in complex analytic case or in neutral signature) and pairs of commuting parameter depending vector fields $X_1(\lambda),X_2(\lambda)$ which are divergence free with respect to some volume form. More precisely, one side of this correspondence is as follows: once
\begin{equation}\label{intro1}\equ
 X_1(\lambda)=V_1+\la V_2, \quad X_2(\lambda)=V_3+\la V_4, \quad [X_1(\la),X_2(\la)]\equiv 0, \quad \mathcal{L}_{V_i}\om=0
 \end{equation}
 for some volume form $\om$, the metric
 \begin{equation}\label{intro2}\equ
 g=\kappa(V^1\odot V^4-V^2\odot V^3),
 \end{equation}
 where $\kappa:=\om(V_1,V_2,V_3,V_4)$ and $(V^i)$ stands for the coframe dual to the frame $(V_i)$, is self-dual and satisfies the vacuum Einstein equations.

Earlier Plebański showed examples of such vector fields depending of one function $f$ of four variables (the potential) satisfying the celebrated I or II Plebański heavenly PDEs \cite{plebanski}:
$$
f_{13}f_{24}-f_{14}f_{23}=1,\qquad f_{13}+f_{24}+f_{11}f_{22}-f_{12}^2=0.
$$
For instance,  the Mason--Newman vector fields and the volume form corresponding to the I Plebański equation are given by
$$
X_1(\la)=-\d_3+\la(f_{13}\d_2-f_{23}\d_1), \quad X_2(\la)=-\d_4+\la(f_{14}\d_2-f_{24}\d_1), \quad \om=dx^1\wedge\cdots\wedge dx^4.
$$
He also proved that any self-dual Ricci-flat metric is governed by this equation. Other PDEs leading to Mason--Newman vector fields are also known in the literature:
$$
f_{34}+f_{23}f_{14}-f_{13}f_{24}=0
$$
(Husain--Park equation \cite{park}, \cite{husain});
$$
f_{44}+f_{12}f_{34}-f_{13}f_{24}=0
$$
(Grant equation \cite{grant});
 \begin{equation}\label{schief}\equ
(\la_1-\la_2)(\la_3-\la_4)f_{12}f_{34}-(\la_1-\la_3)(\la_2-\la_4)f_{13}f_{24}
+(\la_1-\la_4)(\la_2-\la_3)f_{14}f_{23}=0
\end{equation}
(general heavenly equation, Schief equation \cite{schief}, {\color{black} \cite{doubrovFerapontov}, } \cite{konopSchiefszer}). One of the novelties in the approach
{\color{black}of \cite{konopSchiefszer}} was an observation that among the solutions of the general heavenly equation there is a subclass of functions satisfying the so-called dispersionless Hirota system
 \begin{align*}
f_1f_{23}(\la_2-\la_3)+f_2f_{31}(\la_3-\la_1)+f_3f_{12}(\la_1-\la_2)=0\\
f_1f_{24}(\la_2-\la_4)+f_2f_{41}(\la_4-\la_1)+f_4f_{12}(\la_1-\la_2)=0\\
f_1f_{34}(\la_3-\la_4)+f_3f_{41}(\la_4-\la_1)+f_4f_{13}(\la_1-\la_3)=0\\
f_2f_{34}(\la_3-\la_4)+f_3f_{42}(\la_4-\la_2)+f_4f_{23}(\la_2-\la_3)=0
 \end{align*}
 related to the so-called Veronese webs in 4D. Another novelty was {\color{black}the} use  of common eigenfunctions $\phi_i$ of the vector fields $X_1(\lambda_i),X_2(\lambda_i)$, $i=1,\ldots,4$, where $\la_1,\ldots,\la_4$ are four  arbitrary distinct parameters, as privileged coordinates {\color{black}(}which always leads to (\ref{schief}){\color{black})}.
 It is easy to see that this in particular allows to present the vector fields $X_1(\lambda),X_2(\lambda)$ in the form
 \begin{equation*}\label{NijKKS0}\equ
 X_1(\lambda)=(\widetilde{N}+\la \Id)V_2, \qquad  X_2(\lambda)=(\widetilde{N}+\la \Id)V_4,
 \end{equation*}
 where $\widetilde{N}$ is a Nijenhuis operator given by
 \begin{equation}\label{NijKKS}\equ
 \partial_{\phi_i}\mapsto \la_i\partial_{\phi_i}.
  \end{equation}
 On the other hand, one can understand the family of foliations $\{\F_\la\}$ corresponding to the integrable distribution $D_\la=\langle X_1(\lambda),X_2(\lambda)\rangle$ generated by the vector fields $X_1(\lambda)$ and $X_2(\lambda)$ as a special case of a so-called Kronecker web.

In this paper we take advantage of the web theory understood from the point of view of Nijenhuis operators \cite{pHirota}, \cite{pKronwebs} for interpretation and development  of the theory mentioned. In particular, the main result of this paper is a systematic construction of integrable PDEs in four dimensions, some of which lead to Mason--Newman vector fields $X_1(\la),X_2(\la)$ ({\color{black}i.e. vector fields satisfying (\ref{intro1}) for some volume form $\om$  and linearly independent $V_1,\ldots,V_4$}) and self-dual vacuum Einstein metrics, based on different normal forms of 4-dimensional Nijenhuis operators (see Section \ref{secMain}). An important feature of our approach is simultaneous considering of Hirota-type systems, which are in a sense built in the theory.

Let us explain the main ideas of our approach. By a \emph{web}  on a manifold $M$ we understand a 1-parameter family $\{\F_\la\}_{\la\in\P^1}$ of foliations parametrized by 1-dimensional projective space (recall also that a \emph{classical web} is a finite family of foliations). Following \cite{pHirota}, \cite{pKronwebs} one can define the so-called Kronecker webs $\{\F_\la\}$ introduced in \cite{z1} by means of a partial Nijenhuis operator $N:T\F_\infty\to TM$ which can be treated as the restriction of a Nijenhuis operator $\widetilde{N}:TM\to TM$ to the subbundle $T\F_\infty$ of vectors tangent to the foliation $\F_\infty$. The Kronecker web $\{\F_\la\}$ corresponding to the integrable distribution $D_\la=\langle X_1(\lambda),X_2(\lambda)\rangle$, where $X_1(\la),X_2(\la)$ are  Mason--Newman vector fields, has additional property: locally there exists a potential, i.e. a function $f$ such that the 2-form $\be_\la:=dd_{(\widetilde{N}-\la \Id)^{-1}}f$ annihilates $T\F_\la$; here $d_{\widetilde{N}}$ is the so-called Nijenhuis exterior derivative corresponding to a Nijenhuis operator $\widetilde{N}$. The two-form $\be_\la$ is  of rank 2,  i.e. satisfies the identity\footnote{\label{ftn1} {\color{black} Recall that a two-form $\ga_\la$ quadratically depending on $\la$ is called a Gindikin form \cite{gindikin82}, \cite{gindikin85} if it is closed,  nondegenerate, i.e. $\ga_\la\wedge\ga_\mu\not=0$ for $\la\not=\mu$, and is of rank 2. The 2-form $\be_\la$ is rationally depending on the parameter, but in the case, when $\widetilde{N}$ has constant eigenvalues, is a multiple of a Gindikin form with a constant factor.}}
 \begin{equation}\label{bb}\equ
\be_\la\wedge\be_\la=0,
\end{equation}
equivalent to a second order ``heavenly'' PDE on the function $f$, which is \emph{integrable} in the sense that it  possesses a dispersionless Lax pair \cite{BerjFerKrNov}, i.e.
 vector fields $X_1(\la),X_2(\la)$ depending on the second jet of $f$ and an auxiliary parameter $\la$ such that the
Frobenius integrability condition $[X_1, X_2] \in \langle X_1,X_2\rangle$ holds identically modulo the equation and its differential consequences\footnote{All the vector fields from Lax pairs and ``Lax triples'', see Section \ref{6D}, in this paper are linearly depending on the parameter.}.

For instance, the extension of $N$ to the Nijenhuis operator $\widetilde{N}$ given by formula (\ref{NijKKS}), {\color{black}cf. Remark \ref{rem1}, }leads to the general heavenly equation (\ref{schief}). However, the Nijenhuis operator $\widetilde{N}$ that extends $N$ is highly nonunique and other extensions will lead to other integrable PDEs. This idea was exploited in \cite{pHirota} to obtain integrable deformation of the Hirota dispersionless equation and series of other integrable PDEs corresponding to Veronese webs in dimension 3. Veronese webs are special Kronecker webs $\{\F_\la\}$ with $\codim_M \F_\la=1$. In the four-dimensional situation  Veronese webs are described by a system of PDEs which in our scheme is equvalent to the identity
 \begin{equation}\label{ada}\equ
d\al_\la\wedge\al_\la=0,
\end{equation}
 where $\al_\la:=d_{(\widetilde{N}-\la \Id)^{-1}}f$. Note that (\ref{ada}) implies (\ref{bb}), that is why along with a ``heavenly'' equation there always comes the corresponding Hirota-type system and a function satisfying
the latter necessarily satisfies the former. Geometrical picture standing behind is the following: the Kronecker web $\{\F_\la\}$ annihlated by $\be_\la$ is a ``subweb'' of the  Veronese web $\{\V_\la\}$ annihilated by $\al_\la$, i.e. the leaves of the foliation $\F_\la$ are submanifolds in the leaves of the foliation $\V_\la$.

In more details, the content of the paper is the following. In Section \ref{sec1} we recall the notions of a Kronecker web, a partial Nijenhuis operator and Nijenhuis operator and discuss relations between them. In Section \ref{divfree} we introduce the main notion of this paper, divergence free Kronecker webs. In the four-dimensional case they  coincide with the webs generated by  Mason--Newman vector fields $X_{1,2}(\la)$. We prove that a Kronecker web $\{\F_\la\}$ is divergence free if and only if its tangent spaces are annihilated by a polynomially depending on the parameter $\la$ closed differential $k$-form $\be_\la$ for suitable $k$ and suitable degree in $\la$ (the ``if'' part is proven only in 4D case and for the Kronecker webs of the 3-web type, see Definition \ref{3webdefi}).

Section \ref{secPot} is devoted to the proof of the existence of the potential, i.e. a function which satisfies the corresponding integrable ``heavenly'' equation. We prove this existence for any divergence free Kronecker web in 4D under the assumptions that the partial Nijenhuis operator corresponding to this web can be extended to a Nijenhuis operator which is cyclic {\color{black} and the trace of the operator is constant along $\F_\infty$} {\color{black}(see Corollary \ref{corrop})}. To this end we use the ``$\d\bar{\d}$-lemma'' for  Nijenhuis operators \cite{turielDDN},\cite{bolsKonMat}. In particular, our result generalizes the corresponding result from \cite{konopSchiefszer}, where the existence of the potential was proved directly (i.e. in fact  the ``$\d\bar{\d}$-lemma'' for the Nijenhuis operator (\ref{NijKKS}) was recovered). In both Sections \ref{divfree} and \ref{secPot} we use the Koszul operator on multivector fields dual to the exterior derivative \cite{koszul}.

We start Section \ref{secMain} by proving that, given any function $f$ and an invertible Nijenhuis operator $\widetilde{N}$, conditions (\ref{bb}) and (\ref{ada}) are in fact equivalent to $\be_\infty\wedge \be_\infty=0$ and $d\al_0\wedge\al_0$ respectively, i.e. (\ref{bb}) and (\ref{ada}) are indeed equivalent to a PDE and a system of PDEs on $f$ that are independent of the parameter $\la$. These heavenly PDE and a  Hirota-type system of PDEs are automatically integrable in the sense that they are equivalent to the commutation relations of depending on $\la$ vector fields spanning the corresponding integrable distribution $\ker\be_\la$ or $\ker\al_\la$. We then list the heavenly PDEs, the pairs of parameter depending vector fields forming an involutive system, and the Hirota-type systems corresponding to different pairwise nonequivalent local normal forms of Nijenhuis operators in 4D. For the cases with constant eigenvalues we also calculate volume forms $\om$, with respect to which the Mason--Newman vector fields are divergence-free, and the corresponding metrics obtained by formula (\ref{intro2}).

Section \ref{hyperCR} is devoted to discussion of the {\color{black}4D part of the so-called universal hierarchy} and the corresponding heavenly equation. {\color{black} The universal hierarchy of Mart\'{i}nez Alonso and Shabat \cite{martinezAlShabat} was considered in \cite{ferapontovKruglikov}, where its  ``similarity with the Veronese web hierarchy'' was observed. We show that this similarity is not accidental, since the first three equations of the universal hierarchy are equivalent to the Frobenius integrability condition for a special kind of the parameter depending one-form $\al_\la$ describing a Veronese web in 4D (see Theorem \ref{theoHCR}, Corollary \ref{coroHCR}, and Remark \ref{remHCR}). {\color{black}This approach, which can be easily extended to higher dimensions, was originated in \cite{rigal}}. It is natural to ask, whether the {\color{black}4D part of the universal hierarchy} generates some new heavely  PDE. The answer is negative: in Theorem \ref{theo2HCR} we show that in fact we get $\mathbf{E_1}$-case obtained in Section \ref{secMain}.}

Section \ref{sexa} is aimed to show that, on one hand, there can be solutions of heavenly equations corresponding to normal forms of Nijenhuis operators with nonconstant eigenvalues leading to Mason--Newman vector fields (see Examples \ref{exaI}, \ref{exaII} {\color{black} which concern $\mathbf{A_0}$-case}) but, on the other hand, for a general solution {\color{black} such vector fields do not exist. We end this section with a discussion of the $\mathbf{G_0}$-case, which also corresponds to nonoconstant eigenvalues, in the context of existence of Mason--Newman vector fields.}

In Section \ref{6D} we show that the scheme of Section \ref{secMain} works, i.e. produces integrable PDEs or systems of PDEs on $f$, also in higher dimensions, where by integrability we mean existing of a system of vector fields depending on the second jet of $f$ and on an auxiliary  parameter, commutation relations of which are implied by the original system or single PDE and their differential consequences. We present two examples in dimension 5 coming from condition (\ref{bb}). In this case we get a system of 5 equations and a system of 3 commuting parameter depending vector fields.  Finally we present an example in dimension 6, where we use another condition providing nontriviality of the distribution $\ker\be_\la$, namely vanishing of the pfaffian of the matrix of the 2-form $\be_\la$. This gives a single PDE and a system of 2 commuting parameter depending vector fields.

The concluding Section \ref{concl} is devoted to discussing our results in the context of existing literature, in particular of papers \cite{BerjFerKrNov}, \cite{konopSchiefszer}, and pointing out possible further research.

In Appendix we list the local normal forms of Nijenhuis operators which are used in Section \ref{secMain}.

All objects in this paper are either real or complex analytic.

\section{Kronecker webs and Nijenhuis operators}
\label{sec1}

To introduce Kronecker webs we need the following definition \cite{pHirota}, \cite{pKronwebs}.

\abz\label{PNOdefi}
\begin{defi}\rm
Let  $\F$ be a foliation on a manifold $M$.
We say that a bundle morphism $N\colon T\F\to TM$  is a  {\em   partial Nijenhuis operator} (PNO for short) on $M$ if the following two conditions hold:
\begin{enumerate}
\item[(i)] $[X,Y]_N:=[NX,Y]+[X,NY]-N[X,Y]\in \G(T\F)$ for any $X,Y\in \G(T\F)$ (here $[,]$ stands for the commutator of vector fields on $M$ {\color{black} and $\G(T\F)$ for the space of vector fields tangent to the foliation $\F$});
\item[(ii)] $T_N(X,Y):=[NX,NY]-N[X,Y]_N=0$ for any $X,Y\in\G(T\F)$ (it follows from condition (i) that the second term is correctly defined).
\end{enumerate}
\end{defi}

It is easy to see that, given a PNO $N\colon T\F\to TM$, the generalized distribution $D_\la:=(N-\la I)T\F$ (assume for simplicity that $\dim D_\la(x)$ is constant in $x\in M$, i.e. $D_\la$ is a distribution) is integrable for any $\la\in \overline\K$, i.e. de facto we get a web $\{\F_\la\}$, $T\F_\la=D_\la$; here $I:T\F\to TM$ stands for the canonical inclusion, $\overline\K=\K\cup\{\infty\}=\K\P^1$, $\K$ is a ground field $\R$ or $\C$, {\color{black} $D_\infty:=T\F$}.

In the particular case when $T\F=TM$, we recover classical definition of a Nijenhuis operator. On the other hand, we have the following lemma.

\abz\label{llll}
\begin{lemm} (\cite[Lemma 2.5]{pHirota}) Let $\F$ be a foliation on $M$ and let  $N\colon TM\to TM$ be a Nijenhuis   operator such that for some $\la_0\in \K$ the following two conditions hold: (1) the distribution $B:=(N+\la_0\Id_{TM})T\F$ is integrable; (2) $(N+\la_0\Id_{TM})^{-1}(B)=T\F$. Then  $N|_{T\F}$ is a PNO.
\end{lemm}

We see that   a PNO can have an extension to a Nijenhuis operator. It is important to understand that such an extension can be nonunique.

The notion of a Kronecker web \cite{z1} now can be introduced as follows (cf. \cite[Rem. 6.5]{pKronwebs})

\abz\label{krWEBdefi}
\begin{defi}\rm
A collection of foliations $\{\F_\la\}_{\la\in\K\P^1}$ on a manifold $M$ is a \emph{Kronecker web}, if there exists a PNO $N:T\F_\infty\to TM$ such that
\begin{enumerate}\item the morphism $N-\la I:T\F_\infty\to TM$ at any point of the base manifold is injective for any $\la\in\C$ (in the real category we should pass to the complexification of the tangent space);
\item $(N-\la I)T\F_\infty=T\F_\la$ for any $\la\in\K$.
\end{enumerate}
The number $\dim(T_x\F_\la)=\dim(N-\la I)T_x\F_\infty$ is called the \emph{rank} of the Kronecker web.
\end{defi}

The terminology is motivated by the fact that condition 1 {\color{black} implies} the absence of Jordan blocks in the Jordan--Kronecker decomposition of the pair of operators $N_x,I_x:T_x\F_\infty\to T_xM$: {\color{black} the components of this decomposition are the so-called increasing Kronecker blocks, i.e. there exist bases in the restricted domains and codomains of $N_x,I_x$ corresponding to the components, in which these operators are of the form}
\begin{equation}\equ\label{kronform}
\left[
    \begin{array}{cccc}
      0 & & &  \\
      1 & 0 &  &\\
            & \ddots&\ddots& \\
       & & 1 & 0 \\
       & & & 1 \\
    \end{array}
  \right],
  \left[
    \begin{array}{cccc}
      1 & & &  \\
      0 & 1 &  &\\
      & \ddots&\ddots& \\
       & & 0 & 1 \\
       & & & 0 \\
    \end{array}
  \right]
\end{equation}
(these are $(n+1)\times n$-matrices and we say that this \emph{block is of dimension }$n$).
 In the case when at any $x\in M$ we have a sole Kronecker block  the corresponding Kronecker webs are called \emph{Veronese webs} \cite{gz2} (the annihilating 1-form $\al_\la$, $\al_\la(T\F_\la)=0$, in this case generates a Veronese curve in $\P(T^*_xM)$).

Recall that a \emph{classical 3-web} on a manifold $M^{2n}$ is a triple of foliations $\F_1,\F_2,\F_3$ of rank $n$ in general position. There is a canonical Chern affine connection on $M^{2n}$ associated to such a web \cite{chern}. P.~ Nagy \cite{Nagy} showed that in the case, when this connection is torsionless, there exists a web $\{\F_\la\}$ with a special property such that $\F_{\la_i}=\F_i$, $i=1,2,3$, for some $\la_i$. In fact this property means that the web $\{\F_\la\}$ is Kronecker and this motivates the following definition.

\abz\label{3webdefi}
\begin{defi}\rm
A Kronecker web $\{\F_\la\}$,  $T\F_\la=(N-\la I)T\F_\infty$, of rank $n$ on a manifold $M$ of dimension $2n$ is of  \emph{3-web type} if at any $x\in M$ the Jordan--Kronecker decomposition of the pair of operators $N_x,I_x:T_x\F_\infty\to T_xM$ consists of $n$ 1-dimensional Kronecker blocks.
\end{defi}

The last condition can be rephrased as follows: there exist linearly independent vector fields $X_1,\ldots,X_n$, $Y_1,\ldots,Y_n$ on $M$ such that the corresponding PNO is given by $N:X_i\mapsto Y_i$, $i=1,\ldots,n$. The corresponding Konecker web $\{\F_\la\}$ is given by $T\F_\la=\langle Y_1-\la X_1,\ldots,Y_n-\la X_n\rangle$.

\section{Divergence free Kronecker webs}
\label{divfree}

\abz\label{divfreedefi}
\begin{defi}\rm
Let $\om$ be a volume form on $M$. A Kronecker web $\{\F_\la\}_{\la\in\P^1}$,  $T\F_\la=(N-\la I)T\F_\infty$, of rank $n$ on a manifold $M$ is said to be  \emph{divergence free}\footnote{Note that the term ``divergence free webs'' is known in the literature and is used in the context of classical webs \cite{taba}, \cite{domZub}. It means a classical web equipped with a volume form, hence it is not in contradiction with our terminology.} with respect to $\om$ if there exist  vector fields $X_1,\ldots,X_n$    such that for any $\la\in\P^1$ the vector fields $X_1(\la),\ldots,X_n(\la)$,\, $X_i(\la):=(N-\la I)(X_i)$,
\begin{enumerate}\item  are linearly independent and span $T\F_\la$;
\item  pairwise commute: $[X_i(\la),X_j(\la)]=0$;
\item  are divergence free with respect to $\om$: $d(X_i(\la)\iprod\om)=0$.
\end{enumerate}
\end{defi}


\abz\label{th1}
\begin{theo}
Let a Kronecker web $\{\F_\la\}_{\la\in\P^1}$,  $T\F_\la=(N-\la I)T\F_\infty$, of rank $n$ on a manifold $M$  of dimension $m$ be divergence free with respect to a volume form  $\om$ and let $X_1,\ldots,X_n$   be as above. Then there exists a $\la$-depending closed $(m-n)$-form $\beta_\la$ annihilating $T\F_\la$ (i.e. $X_i(\la)\iprod \beta_\la=0$ for any $\la$). It is given by
$$
\beta_\la:=X_1(\la)\wedge \cdots\wedge X_n(\la)\iprod \om
$$
and is a polynomial in $\la$ of degree $n$.
\end{theo}

\noindent{\sc Proof} The fact that $\be_\la$ annihilates $T\F_\la$ and is a polynomial of degree $n$ is obvious. To prove the closedness of $\beta_\la$ we will exploit the Koszul operator \cite{koszul}
$$
D=\Phi^{-1}\circ d\circ \Phi:\chi^k(M)\to\chi^{k-1}(M),
$$
where $\chi^k(M)$ stands for the space of $k$-vector fields on $M$ and $\Phi$ is the isomorphism of  $\chi^k(M)$ with the space of $(m-k)$-forms given by the contraction with\footnote{\label{ftn}More precisely,  we  {\color{black} define $\Phi$ by} $\Phi(U)=i(U)\om$, where $U\in \chi^k(M)$ and the interior product $i(U)\om$ is  {\color{black} given} by $\langle i(U)\om,V\rangle=\langle \om, U\wedge V\rangle$, $V\in\chi^{m-k}(M)$, $\langle ,\rangle$ being the natural pairing between forms and multivectors.} $\om$. It is related with the Schouten bracket of multivector fields by the formula
\begin{equation}\label{DKos}\equ
[U,V]=(-1)^k (D(U\wedge V)-D(U)\wedge V-(-1)^k U\wedge D(V)), \qquad  U\in\chi^k(M), \quad V\in\chi^l(M).
\end{equation}
Since $D(X_i(\la))=0$ and the vector fields pairwise commute by induction we get that $D(X_1(\la)\wedge \cdots\wedge X_n(\la))=0$. \qed

\medskip

In the case of Kronecker webs of 3-web type in 4-dimensional case we can prove also the converse statement.

\abz\label{th1a}
\begin{theo}
Let a Kronecker web $\{\F_\la\}_{\la\in\P^1}$,  $T\F_\la=(N-\la I)T\F_\infty$, of rank $2$ on a manifold $M^{4}$ be of 3-web type. Assume that there exists a $\la$-depending closed $2$-form $\beta_\la$ annihilating $T\F_\la$ such that the dependence of $\be_\la$ on $\la$ is polynomial of second order. Then $\{\F_\la\}$ is divergence free with respect to some volume form $\om$.
\end{theo}

\noindent{\sc Proof}  Choose coordinate system $(x_1,x_2,y_1,y_2)$ on $M$ such that $\langle \partial_{x_1},\partial_{x_2}\rangle=T\F_\infty$ and $\langle \partial_{y_1},\partial_{y_2}\rangle=T\F_0$  and put $X_1:=\partial_{x_1},X_2:=\partial_{x_2}$. Then
$$
NX_1=X_{11}\partial_{y_1}+X_{12}\partial_{y_2}, \qquad NX_2=X_{21}\partial_{y_1}+X_{22}\partial_{y_2}
$$
for some functions $X_{ij}$ such that $\Delta:=X_{11}X_{22}-X_{12}X_{21}$ {\color{black} is nonvanishing}. Since the distribution $\langle X_1(\la),X_2(\la) \rangle$, $X_i(\la):=(N-\la I)X_i$, is integrable for any $\la$, it follows from the form of the vector fields $X_1(\la),X_2(\la)$ that they commute.

One checks directly that the 2-form
$$
\be_\la=a( \Delta dx_1\wedge dx_2-\la(X_{12} dx_1\wedge dy_1-X_{11}dx_1\wedge dy_2+X_{22} dx_2\wedge dy_1-X_{21}dx_2\wedge dy_2)+ \la^2dy_1\wedge dy_2),
$$
where $a$ is some function independent of $\la$, annihilates $X_1(\la)$ and $X_2(\la)$. Note that any 2-form of second order in $\la$ annihilating $X_1(\la)$ and $X_2(\la)$ is of this form.  Now put
\begin{equation}\label{omm}\equ
\om=\frac1{\la^2}\be_\la\wedge dx_1 \wedge dx_2=a dy_1 \wedge dy_2 \wedge dx_1 \wedge dx_2.
\end{equation}
 Then $d(X_1(\la)\iprod \om)=-\frac1{\la}d(\be_\la\wedge dx_2)=0$ and $d(X_2(\la)\iprod \om)=\frac1{\la}d(\be_\la\wedge dx_1)=0$. \qed

\medskip

{\color{black}Combining Theorems \ref{th1} and \ref{th1a} we get the following}

\abz\label{rem02}
\begin{coro}
A Kronecker web $\{\F_\la\}_{\la\in\P^1}$ on a $4$-dimensional manifold $M$  of 3-web type is divergence free with respect to some volume form if and only if  there exists a $\la$-depending closed $2$-form $\beta_\la$ annihilating $T\F_\la$ such that the dependence of $\be_\la$ on $\la$ is polynomial of second order.
\end{coro}

{\color{black}
\abz\label{rema003}
\begin{rema}\rm Existence of a closed quadratically depending on the parameter 2-form $\be_\la$ annihilating the Mason--Newman vector fields is not new in the literature. For instance such forms is an essential tool in studying vacuum solutions of the Einstein equation in \cite{gindikin82}, \cite{gindikin85} {\color{black} cf. footnote \ref{ftn1} on page \pageref{ftn1}}. Also, such a form induces a fiber-closed 2-form on the twistor fibration $\mathcal{PT}\to\P^1$ corresponding to a SD vacuum Einstein metric and plays a fundamental role in the reconstruction of the metric from the twistor space \cite{penrose}, \cite[Th. 10.5.5]{dunajskiBook}.
\end{rema}
}

\section{Nijenhuis exterior derivative and existence of the potential}
\label{secPot}

We will exploit the following facts about the so-called Nijenhuis differential $d_N$ on differential forms \cite{mks} given by the formula $d_N=[i_N,d]$, where
$$
(i_N\al)(X_1,X_2,\ldots,X_k)=\al(NX_1,X_2\ldots,X_k)+\al(X_1,NX_2\ldots,X_k)+\cdots+\al(X_1,X_2\ldots,N
X_k);
$$
here $N:TM\to TM$ is a Nijenhuis operator, $\al$ is a differential $k$-form on $M$ (by definition $i_N\al=0$ for a $0$-form $\al$).

{\color{black} The following theorem is a crucial ingredient in our considerations below concerning the existence of the potential.}

\abz\label{th2}
\begin{theo}(\cite{turielDDN},\cite{bolsKonMat})
Let $N:TM\to TM$ be a Nijenhuis operator such that the operator $N_x:T_xM\to T_xM$ is cyclic for any $x\in M$, i.e. the Frobenius normal form of this operator consists of one cyclic block. Let $\be$ be a differential 2-form which is $d$- and $d_N$-closed. Then locally there exists a function $\T$ such that
$$
dd_N\T=\be.
$$
\end{theo}

Theorem \ref{th3} below generalizes the Koszul formula (\ref{DKos}) to the case of the Nijenhuis differential (we consider only the case of vector fields). Let $N:TM\to TM$ be a Nijenhuis operator and  $\om$  a volume form on $M$.
Introduce the Koszul--Nijenhuis operator by
$$
D_N=\Phi^{-1}\circ d_N\circ \Phi:\chi^k(M)\to\chi^{k-1}(M),
$$
where as previously $\chi^k(M)$ stands for the space of $k$-vector fields on $M$ and $\Phi$ is the isomorphism of  $\chi^k(M)$ with the space of $(m-k)$-forms given by the contraction with $\om$.
\abz\label{th3}
\begin{theo}
If $U,V$ are vector fields on $M$, then
$$
[U,V]_N=-D_N(U\wedge V)+D_N(U)\wedge V-U\wedge D_N(V),
$$
where $[,]_N$ is the Lie bracket on the space of vector fields on $M$ given by the formula from Definition \ref{PNOdefi}.
\end{theo}

\abz\label{lem12}\equ
\begin{lemm}
Let $I_N:=\Phi^{-1}\circ i_N\circ\Phi:\chi^k(M)\to \chi^k(M)$. Then
\begin{equation}\label{eIN}
I_N=-j_{N}+\Tr(N)\Id_{\chi^k(M)},
\end{equation}
where
$$
j_N(X_1\wedge X_2\wedge\cdots \wedge X_k):=NX_1\wedge X_2\wedge\cdots \wedge X_k+X_1\wedge NX_2\wedge\cdots\wedge X_k+\cdots+X_1\wedge X_2\wedge\cdots \wedge NX_k.
$$
Moreover,
\begin{equation}\label{eDN}\equ
D_N(U)=-[j_N,D](U)-[\Tr(N),U], \qquad U\in\chi^k(M),
\end{equation}
where $[,]$ in the last term stands for the Schouten bracket.
\end{lemm}

\noindent{\sc Proof} Let $U\in\chi^k(M)$. Without loss of generality we can assume that locally $U=X_1\wedge \cdots \wedge X_k$ for some vector fields $X_i$. Let $X_1,\ldots, X_k,X_{k+1},\ldots, X_m$ be a local frame on $M$ and $X^1,\ldots, X^m$ the dual coframe. Then $\om=fX^1\wedge\cdots\wedge X^m$ for some function $f$ and $\Phi(U)=fX^{k+1}\wedge\cdots\wedge X^m$.

Let
$$
NX_i=N_i^jX_j+N_i^\be X_\be, NX_\al=N_\al^jX_j+N_\al^\be X_\be,
$$
where we use the summation convention and we denote the indices from $\{1,\ldots,k\}$ by Latin letters and that from $\{k+1,\ldots,m\}$ by Greek letters. Then also
$$
N^tX^\al=N_j^\al X^j+N_\be^\al X^\be,
$$
where we denote by $N^t$ the transposed to $N$ operator, and
\begin{align*}
j_N(U)=N_1^\be X_\be\wedge X_2\wedge\cdots \wedge X_k+X_1\wedge N_2^\be X_\be\wedge\cdots\wedge X_k+\cdots+X_1\wedge X_2\wedge\cdots \wedge N_k^\be X_\be+\\
(N_1^1+N_2^2+\cdots+N_k^k)U.
\end{align*}
Thus  (cf. footnote \ref{ftn} on page \pageref{ftn})
\begin{align*}
\Phi\circ j_N(U)=f[(-1)^{k-1}\sum_\be (-1)^{\be-1}N_1^\be X^1\wedge X^{k+1}\wedge\cdots\wedge \check{X^\be}\wedge\cdots \wedge X^m+\\
(-1)^{k-2}\sum_\be (-1)^{\be-2}N_2^\be X^2\wedge X^{k+1}\wedge\cdots\wedge \check{X^\be}\wedge\cdots \wedge X^m+\cdots+\\
(-1)^{k-k}\sum_\be (-1)^{\be-k}N_k^\be X^k\wedge X^{k+1}\wedge\cdots\wedge \check{X^\be}\wedge\cdots \wedge X^m+\\
(N_1^1+N_2^2+\cdots+N_k^k)X^{k+1}\wedge\cdots\wedge X^m]=\\
f[\sum_i(-1)^{k-i}\sum_\be (-1)^{\be-i}(-1)^{\be-k-1} X^{k+1}\wedge\cdots\wedge X^{\be-1}\wedge N_i^\be X^i\wedge X^{\be+1}\wedge\cdots \wedge X^m+\\
(N_1^1+N_2^2+\cdots+N_k^k)X^{k+1}\wedge\cdots\wedge X^m]=\\
f[\sum_i\sum_\be (-1)^{2\be-2i-1} X^{k+1}\wedge\cdots\wedge X^{\be-1}\wedge N_i^\be X^i\wedge X^{\be+1}\wedge\cdots \wedge X^m+\\
(N_1^1+N_2^2+\cdots+N_k^k)X^{k+1}\wedge\cdots\wedge X^m]=\\
f[-\sum_\be X^{k+1}\wedge\cdots\wedge X^{\be-1}\wedge N^tX^\be\wedge X^{\be+1}\wedge\cdots \wedge X^m+\\
(N_{k+1}^{k+1}+N_{k+2}^{k+2}+\cdots+N_m^m)X^{k+1}\wedge\cdots\wedge X^m+(N_1^1+N_2^2+\cdots+N_k^k)X^{k+1}\wedge\cdots\wedge X^m]=\\
-f[N^tX^{k+1}\wedge X^{k+2}\wedge\cdots\wedge X^m+X^{k+1}\wedge N^t X^{k+2}\wedge\cdots\wedge X^m+X^{k+1}\wedge X^{k+2}\wedge\cdots\wedge N^tX^m]+\\
f[\Tr(N)X^{k+1}\wedge\cdots\wedge X^m]=
(-i_{N}+\Tr(N)\Id_{\chi^k(M)})\circ \Phi(U),
\end{align*}
where $\check{X^\be}$ as usual denotes the omitted component. Formula (\ref{eIN}) is proved. To prove formula (\ref{eDN})
observe that by definition and by (\ref{eIN}) we have $D_N=\Phi^{-1}\circ d_N\circ\Phi=\Phi^{-1}\circ [i_N,d]\circ\Phi= [I_N,D]=-[j_N,D]+[\Tr(N)\Id_{\chi^k(M)},D]$. Finally it follows from (\ref{DKos}) that $[g\Id_{\chi^k(M)},D](U)=-[g,U]$ for $U\in\chi^k(M)$ and any function $g$. \qed

\medskip

\noindent{\sc Proof of Theorem \ref{th3}} Put $g:=\Tr(N)$. Note that for any $U,V\in\chi^1(M)$ we have
\begin{align*}
[g,U\wedge V]-[g,U]\wedge V+U\wedge [g,V]=0
\end{align*}
due to the standard properties of the Schouten bracket. Thus it follows from (\ref{eDN}) that
\begin{align*}
-D_N(U\wedge V)+D_N(U)\wedge V-U\wedge D_N(V)=\\
[j_N,D](U\wedge V)-[j_N,D](U)\wedge V+U\wedge [j_N,D](V)=\\
ND(U\wedge V)-D(NU\wedge V+U\wedge NV)+D(NU) V-UD(NV)= \\
ND(U\wedge V)+[NU,V]+[U,NV]+NU\cdot D(V)-D(U)\cdot NV=\\
[NU,V]+[U,NV]-N[U,V]=[U,V]_N;
\end{align*}
here we used (\ref{DKos}) and the fact that $j_N(f)=0$ for any function $f$ by definition. \qed

\medskip

The main results of this section are the following theorem {\color{black}and its corollary below.}

\abz\label{th4}
\begin{theo}
Let a Kronecker web $\{\F_\la\}_{\la\in\P^1}$,  $T\F_\la=(N-\la I)T\F_\infty$, of rank $2$ on a manifold $M$ of dimension $4$ be divergence free with respect to a volume form  $\om$ and let $X_1,X_2$   be vector fields from  Definition \ref{divfreedefi}, i.e. such that for any $\la\in\P^1$ the vector fields $X_1(\la),X_2(\la)$,\, $X_i(\la):=(N-\la I)(X_i)$,
\begin{enumerate}\item  are linearly independent and span $T\F_\la$;
\item commute: $[X_1(\la),X_2(\la)]=0$;
\item  are divergence free with respect to $\om$: $d(X_i(\la)\iprod\om)=0$.
\end{enumerate}
Assume that  $\widetilde{N}:TM\to TM$ is a Nijenhuis operator such that $\widetilde{ N}|_{T\F_\infty}=N$.

Then the 2-form
\begin{equation}\label{forbe}\equ
\be:=i(X_1\wedge X_2)\om
\end{equation}
(see footnote \ref{ftn} on page \pageref{ftn})
\begin{enumerate}\item[(i)]
is closed with respect to  $d$;
\item[(ii)]
is closed with respect to $d_{\widetilde{N}}$ if and only if
the function $\Tr(\widetilde{N})$ is constant along $\F_\infty$, i.e.
\begin{equation}\label{tra}\equ
X_i\iprod d(\Tr(\widetilde{ N}))=0, \qquad i=1,2.
\end{equation}
\end{enumerate}

\end{theo}

\noindent{\sc Proof} $d$-closedness of $\be$ is equivalent to $D$-closedness of the bivector $X_1\wedge X_2$ and follows from Defintion \ref{divfreedefi}(2,3) and the Koszul formula (\ref{DKos}). Analogously, $d_{\widetilde{ N}}$-closedness of $\be$ is equivalent to $D_{\widetilde{ N}}$-closedness of the bivector $X_1\wedge X_2$. From Defintion \ref{divfreedefi}(2) we conclude that $[X_1,X_2]=0,[NX_1,NX_2]=0,[NX_1,X_2]+[X_1,NX_2]=0$, hence $[X_1,X_2]_{\widetilde{N}}=0$. Thus by Theorem \ref{th3} and by formula (\ref{eDN}) we have
\begin{align*}
0&=-D_{\widetilde{N}}(X_1\wedge X_2)+D_{\widetilde{N}}(X_1)\wedge X_2-X_1\wedge D_{\widetilde{N}}(X_2)\\
&=-D_{\widetilde{N}}(X_1\wedge X_2)-[j_{\widetilde{N}},D](X_1)\wedge X_2+X_1\wedge [j_{\widetilde{N}},D](X_2)
-[\Tr(\widetilde{N}),X_1]\wedge X_2+X_1\wedge [\Tr(\widetilde{N}),X_2)]\\
&=-D_{\widetilde{N}}(X_1\wedge X_2)+D(NX_1)\wedge X_2-X_1\wedge D(NX_2)
-[\Tr(\widetilde{N}),X_1] X_2+ [\Tr(\widetilde{N}),X_2]X_1\\
&=-D_{\widetilde{N}}(X_1\wedge X_2)-[\Tr(\widetilde{N}),X_1] X_2+ [\Tr(\widetilde{N}),X_2]X_1;
\end{align*}
here we used the fact that $D(NX_i)=0$ which follows from assumption 3 (for $\la=0$). Hence,  we have
$$
D_{\widetilde{N}}(X_1\wedge X_2)=-[\Tr(\widetilde{N}),X_1] X_2+ [\Tr(\widetilde{N}),X_2]X_1
$$
and in view of linear independence of $X_1,X_2$ we get \emph{(ii)}.
\qed

\abz\label{corrop}
\begin{coro} Retaining the assumptions of Theorem \ref{th4} assume additionally
that  $\widetilde{N}_x:T_xM\to T_xM$ is cyclic for any $x\in M$, {\color{black}i.e. the Frobenius normal form of this operator consists of one cyclic block}.

Then,
if (\ref{tra}) holds,  locally there exists a function $\Theta$ such that
$$
\be=dd_{\widetilde{ N}}\Theta.
$$
Moreover, the 2-form $\beta_\la=i(X_1(\la)\wedge X_2(\la))\om$ annihilating the vector fields $X_1(\la),X_2(\la)$ is proportional to the form
$$
dd_{\widetilde{ N}-\la \Id_{TM}}\Theta.
$$
\end{coro}

\noindent{\sc Proof}
If (\ref{tra}) holds, we have  $d$- and $d_{\widetilde{N}}$-closedness of the form $\be$ and Theorem \ref{th2} can be used. The second statement of the corollary follows from Lemma \ref{th56} below.\qed

\abz\label{rem2}
\begin{rema}\rm We conjecture that a generalization of Theorem \ref{th4} is valid for any divergence free Kronecker web of corank 2.
\end{rema}

\abz\label{rem1}
\begin{rema}\rm In \cite{konopSchiefszer} the existence of the potential $\Theta$ is proved directly for the extension $\widetilde{N}$ of the PNO $N$ given by
$$
\widetilde{N}:\partial_{\phi_j}\mapsto \la_j\partial_{\phi_j}, \qquad j=1,\ldots,4,
$$
where $\la_1,\ldots,\la_4$ are arbitrary pairwise distinct values of the parameter $\la$ and $\phi_j$ are four local independent functions satisfying
\begin{equation}\label{eX12}\equ
X_1(\la_j)\phi_j=((N-\la_j I)X_1)\phi_j=0, \qquad X_2(\la_j)\phi_j=((N-\la_j I)X_2)\phi_j=0.
\end{equation}
The fact that $\widetilde{N}$ indeed is an extension of the PNO $N:X_i\to Y_i:=NX_i$ is obvious: if $X_i=X_i^j\partial_{\phi_j}$, $Y_i=Y_i^j\partial_{\phi_j}$, $i=1,2$, then (\ref{eX12}) implies $Y_i^j=\la_jX_i^j$ (no summation) and $\widetilde{N}X_i=Y_i$.
\end{rema}

{\color{black}
\abz\label{rem42}
\begin{rema}\rm Consider the 2-form $\be_\la=\be_0+\la\be_1+\la^2\be_2=i(X_1(\la)\wedge X_2(\la))\om$ annihilating the vector fields $X_1(\la),X_2(\la)$.
The 2-form $\be$ given by (\ref{forbe}), which is $dd_{\widetilde{N}}$-closed and is used for the proof of the existence of the potential, is obviously equal to $\be_2$. Notice that
the existence of the potential within the framework of the Plebański I equation is related with the complexification of the ``usual'' $\partial\bar{\partial}$-lemma from K\"{a}hler geometry applied to the form $\be_1$, see \cite[Ch. 13.4]{masonWoodhouse}.
\end{rema}
}

\section{Divergence free Kronecker webs and heavenly PDEs corresponding to different normal forms of Nijenhuis operators in 4D. Relation with  Hirota-type dispersionless systems.}
\label{secMain}

\abz\label{lem56}
\begin{defi}\rm
By a \emph{kernel} of a differential $k$-from $\be$ we understand the  distribution\footnote{Here the anlyticity of the objects is crucial: the set of points in $x\in M$, where rank of $\be(x)$ is nonmaximal, is analytic and by continuity $\rk \ker\be(x)=const$.} $\ker\be$ on $M$ generated by the local vector fields $X$ such that $X\iprod \be=0$.
\end{defi}

Let $N:TM\to TM$ be a Nijenhuis operator, $f$ a function on $M$. Put  $N_\la:=N-\la \Id_{TM}$ and assume that there exist an open set $U \subset M$ and at most finite set $E:=\{\la_1,\ldots,\la_k\}\subset \P^1$ , $\la_i\not=0$ ($E$ may be empty) such that $N_\la$ is invertible for any $\la\in\P^1\setminus E$ on an open dense subset in  $U$. Put
\begin{align*}
 \al_\la:&=d_{(N_\la)^{-1}}f=((N_\la)^{-1})^tdf,   &\al_\infty&:=df\\
 \be_\la:&=d\al_\la=d((N_\la)^{-1})^tdf,  &\be_\infty&:=dd_Nf=dN^tdf,
\end{align*}
where as previously $(\cdot)^t$ stands for the transposed operator.
Consider the following conditions over $U$:
\begin{equation}\label{cond1}\equ
\forall\la\not=\la_i\ \be_\la\wedge  \be_\la=0,
\end{equation}
\begin{equation}\label{cond2}\equ
 \be_\infty\wedge  \be_\infty=0,
\end{equation}
\begin{equation}\label{cond3}\equ
\forall\la\not=\la_i\ \al_\la\wedge  \be_\la=0,
\end{equation}
\begin{equation}\label{cond4}\equ
 \al_0\wedge  \be_0=0.
\end{equation}

\abz\label{th561}
\begin{theo}
Conditions (\ref{cond1}) and (\ref{cond2}), as well as (\ref{cond3}) and (\ref{cond4}) are equivalent.
\end{theo}
\noindent{\sc Proof} Equivalence (\ref{cond1}) $\Longleftrightarrow$ (\ref{cond2})  follows from Lemma \ref{th56} below. Indeed, condition (\ref{cond1}), respectively (\ref{cond2}), means that rank of the 2-form $\be_\la$, respectively $\be_\infty$, at generic point is{\footnote{\color{black}By rank of a 2-form we understand codimension of its kernel.} 2, i.e. rank of the distribution $\ker\be_\la$,  respectively  $\ker\be_\infty$, is $\dim M-2$ (unless $\be_\la$, $\be_\infty$ are nontrivial). By Lemma \ref{th56}   the last two conditions are equivalent.

Implication (\ref{cond3}) $\Longrightarrow$ (\ref{cond4}) is obvious. Finally, to prove the inverse implication write $\F$ for the foliation cut by $f$, i.e. $T\F=\ker df=\ker \al_\infty$, and observe that
condition  (2) of Lemma \ref{llll} holds for $\la_0=0$  due to the invertibility of $N$. Since (\ref{cond4}) means the integrability of the distribution $\ker\al_0=N\ker\al_\infty$, then condition (1) of this lemma also holds and we conclude that $N|_{T\F}$ is a partial Nijenhuis operator and the distribution $(N-\la{\color{black}\Id})T\F$ (which is  equal to $\ker\al_\la$ for any $\la\not=\la_1,\ldots,\la_k$) is integrable for any $\la$, so in fact (\ref{cond3}) holds.  \qed

\medskip

\abz\label{th56}
\begin{lemm}
Put  $D_\la:=\ker\be_\la$, $D_\infty=\ker\be_\infty$.  Then
$$
D_\la=N_\la(D_\infty).
$$
\end{lemm}

\noindent{\sc Proof} Let  $X$ be a vector field. We have to prove the equivalence
$$
X\iprod dN^tdf=0\quad\Longleftrightarrow\quad N_\la X\iprod d((N_\la)^{-1})^tdf=0.
$$
Indeed, let $Y$ be a vector field. Then (since $N_\la$ is also a Nijenhuis operator) we have
\begin{align*}
(N_\la X\iprod d((N_\la)^{-1})^tdf)(N_\la Y)=(d((N_\la)^{-1})^tdf)(N_\la X,N_\la Y)=\\
N_\la Xdf(Y)-N_\la Ydf(X)-df(N_\la^{-1}[N_\la X,N_\la Y])=\\
N_\la Xdf(Y)-N_\la Ydf(X)-df([ X, Y]_{N_\la})=\\
N_\la Xdf(Y)-Ydf(N_\la X)-df([N_\la X,Y])\\
+Xdf(N_\la Y)-N_\la Ydf(X)-df([X, N_\la Y])\\
-Xdf(N_\la Y) +Ydf(N_\la X)+df(N_\la[X,Y])=\\
ddf(N_\la X,Y)+ddf(X,N_\la Y)+d(N_\la ^t df)(X,Y)=\\
(X\iprod d(N_\la ^t df))(Y)=(X\iprod d(N^t df))(Y)-\la (X\iprod d( df))(Y)=\\
(X\iprod d(N^t df))(Y).
\end{align*}
\qed

\medskip

The following list contains: (1) the heavenly PDEs equivalent to condition (\ref{cond2}) corresponding to the local models of Nijenhuis operators in 4D listed in Appendix; (2) the pairs of parameter depending vector fields $X_1(\la),X_2(\la)$, either commuting or forming an involutive system, which span the distribution $\ker\be_\la$, and their commutator; (3) the corresponding Hirota-type systems equivalent to condition (\ref{cond4})\footnote{\label{ft43}Note that only three of four equations (these last coming from the four coefficients of the 3-form $\al_0\wedge\be_0$) in each of these systems are independent.}; (4) the corresponding integrable 1-forms $\al_\la$; (5)
in case of local models of Nijenhuis operators with constant eigenvalues (cf. discussion in Section \ref{sexa}), the corresponding volume forms $\om$, with respect to which $X_{1,2}(\la)$ are divergence free,  and the self-dual Einstein metrics $g$ obtained by formula (\ref{intro2}) within the Mason--Newman formalism.

For short we will exploit only the extreme cases from the list of cyclic normal forms from Appendix: the one when the eigenvalues of $N$ form a functionally independent system and that when all of them are constant, i.e. the cases $\mathbf{A_0,A_4}$, $\mathbf{B_0,B_4}$, $\mathbf{C_0,C_3}$, $\mathbf{D_0,D_3}$, and $\mathbf{E_0,E_1}$.
 Note however that the intermediate cases $\mathbf{A_1,A_2,A_3}$, $\mathbf{B_1,B_2,B_3}$, $\mathbf{C_1, C_2}$, $\mathbf{D_1, D_2}$ also give integrable PDEs and it is possible to explicitly write the parameter depending pairs of vector fields forming an involutive system and the corresponding Hirota systems.

Surprisingly, the following non-cyclic normal forms listed in Appendix also produce integrable PDEs and we include them in the list: $\mathbf{F_0,G_0}$ and $\mathbf{H_0}$.

We checked also some other obvious non-cyclic normal forms\footnote{Note that it is impossible to give an exhausting list of non-cyclic normal forms since such would include all the singularities of Nijenhuis operators, which are not known in general.} and found that they lead to degenerate cases in which the restriction $N|_D$ of the Nijenhuis operator $N$ is not a Kronecker PNO, i.e. the integrable distribution $D_\la=\ker d((N_\la)^{-1})^tdf$ does not correspond to a Kronecker web of 3-web type.

 \begin{enumerate}
\item[$\bf(A_0)$]
$\boxed{(x_1-x_2)(x_3-x_4)f_{12}f_{34}-(x_1-x_3)(x_2-x_4)f_{13}f_{24}+(x_1-x_4)(x_2-x_3)f_{14}f_{23}=0}$
\begin{align}
X_1(\la)=\frac{f_{24}(x_2-x_4)}{f_{12}(x_1-x_2)}(x_1-\la)\partial_{ x_1}-\frac{f_{14}(x_1-x_4)}{f_{12}(x_1-x_2)}(x_2-\la)\partial_{ x_2}+
(x_4-\la)\d_{ x_4}\nonumber \\
 X_2(\la)=\frac{f_{23}(x_2-x_3)}{f_{12}(x_1-x_2)}(x_1-\la)\partial_{ x_1}-\frac{f_{13}(x_1-x_3)}{f_{12}(x_1-x_2)}(x_2-\la)\partial_{ x_2}+
(x_3-\la)\d_{ x_3}\nonumber\\
[X_1(\la),X_2(\la)]=0\nonumber
 \end{align}
  \begin{equation} \equ\label{HirotaA0}
     \begin{split}
f_1f_{23}(x_2-x_3)+f_2f_{31}(x_3-x_1)+f_3f_{12}(x_1-x_2)=0\\
f_1f_{24}(x_2-x_4)+f_2f_{41}(x_4-x_1)+f_4f_{12}(x_1-x_2)=0\\
f_1f_{34}(x_3-x_4)+f_3f_{41}(x_4-x_1)+f_4f_{13}(x_1-x_3)=0\\
f_2f_{34}(x_3-x_4)+f_3f_{42}(x_4-x_2)+f_4f_{23}(x_2-x_3)=0
     \end{split}
 \end{equation}
  \begin{align*}
\al_\la=\frac{f_1}{x_1-\la}d{x_1}+
\frac{f_2}{x_2-\la}d{x_2}+\frac{f_3}{x_3-\la}d{x_3}+
\frac{f_4}{x_4-\la}d{x_4}
 \end{align*}
 \item[$\bf(A_4)$]\label{caseHirota}
$\boxed{(\la_1-\la_2)(\la_3-\la_4)f_{12}f_{34}-(\la_1-\la_3)(\la_2-\la_4)f_{13}f_{24} +(\la_1-\la_4)(\la_2-\la_3)f_{14}f_{23}=0}$
\begin{align*}
X_1(\la)=\frac{f_{24}(\la_2-\la_4)}{f_{12}(\la_1-\la_2)}(\la_1-\la)\partial_{ x_1}-\frac{f_{14}(\la_1-\la_4)}{f_{12}(\la_1-\la_2)}(\la_2-\la)\partial_{ x_2}+
(\la_4-\la)\d_{ x_4}\\
X_2(\la)=\frac{f_{23}(\la_2-\la_3)}{f_{12}(\la_1-\la_2)}(\la_1-\la)\partial_{ x_1}-\frac{f_{13}(\la_1-\la_3)}{f_{12}(\la_1-\la_2)}(\la_2-\la)\partial_{ x_2}+
(\la_3-\la)\d_{ x_3}\\
[X_1(\la),X_2(\la)]=0
 \end{align*}
 \begin{align*}
f_1f_{23}(\la_2-\la_3)+f_2f_{31}(\la_3-\la_1)+f_3f_{12}(\la_1-\la_2)=0\\
f_1f_{24}(\la_2-\la_4)+f_2f_{41}(\la_4-\la_1)+f_4f_{12}(\la_1-\la_2)=0\\
f_1f_{34}(\la_3-\la_4)+f_3f_{41}(\la_4-\la_1)+f_4f_{13}(\la_1-\la_3)=0\\
f_2f_{34}(\la_3-\la_4)+f_3f_{42}(\la_4-\la_2)+f_4f_{23}(\la_2-\la_3)=0
 \end{align*}
  \begin{align*}
\al_\la=\frac{f_1}{\la_1-\la}dx_1+
\frac{f_2}{\la_2-\la}dx_2+\frac{f_3}{\la_3-\la}dx_3+
\frac{f_4}{\la_4-\la}dx_4
 \end{align*}
\begin{align*}
\om:=\frac{(\la_1-\la)(\la_2-\la)}{(\la_1-\la_2)}d\al_\la\wedge dx_3\wedge dx_4=f_{12}dx_1\wedge dx_2\wedge dx_3\wedge dx_4
 \end{align*}
 $$
 g:=\frac1{J}\left[\begin{array}{cccc}
2f_{12}f_{13}f_{14} & f_{12}(f_{14}f_{23}+f_{13}f_{24}) & f_{13}(f_{14}f_{23}+f_{12}f_{34}) & f_{14}(f_{13}f_{24}+f_{12}f_{34})\\
* & 2 f_{12}f_{23}f_{24} & f_{23}(f_{13}f_{24}+f_{12}f_{34}) & f_{24}(f_{14}f_{23}+f_{12}f_{34}) \\
* & * & 2f_{13}f_{23}f_{34} & f_{34}(f_{14}f_{23}+f_{13}f_{24})\\
* & * & * & 2f_{14}f_{24}f_{34}
\end{array}\right],
 $$
 where $J=f_{14}f_{23}-f_{13}f_{24}$.
 \item[$\bf(B_0)$]
$\boxed{(x_1-x_3)(x_1-x_4)(f_{13}f_{24}-f_{23}f_{14})+(x_3-x_4)(f_{23}f_{24}-(f_{22}-f_{2})f_{34})=0}$
\begin{align*}
X_1(\la)=\frac{f_{24}(x_1-x_4)}{f_{22}-f_2}(x_1-\la)\partial_{ x_1}+\frac{f_{24}(x_1-x_4)-(f_{14}(x_1-x_4)+f_{24})(x_1-\la)}{f_{22}-f_2}\partial_{ x_2}+
(x_4-\la)\d_{ x_4}\\
X_2(\la)=\frac{f_{23}(x_1-x_3)}{f_{22}-f_2}(x_1-\la)\partial_{ x_1}+\frac{f_{23}(x_1-x_3)-(f_{13}(x_1-x_3)+f_{23})(x_1-\la)}{f_{22}-f_2}\partial_{ x_2}
+(x_3-\la)\d_{ x_3}\\
[X_1(\la),X_2(\la)]=0
 \end{align*}
 \begin{align*}
-f_2f_{23}+(x_1-x_3)(f_1 f_{23}- f_2f_{13})+f_{3}(f_{22}-f_2)=0\\
-f_2f_{24}+(x_1-x_4)(f_1 f_{24}- f_2f_{14})+f_{4}(f_{22}-f_2)=0\\
(x_3-x_4)(x_1f_1-f_2)f_{34}+x_1(x_1-x_3)f_4f_{13}-x_1(x_1-x_4)f_3f_{14}+x_3f_4f_{23}-x_4f_3f_{24}=0\\
f_4f_{23}(x_1-x_3)-f_3f_{24}(x_1-x_4)+f_2f_{34}(x_3-x_4)=0
 \end{align*}
  \begin{align*}
\al_\la=\frac{(x_1-\la)f_1-f_2}{(x_1-\la)^2}dx_1+
\frac{f_2}{x_1-\la}dx_2+\frac{f_3}{x_3-\la}dx_3+
\frac{f_4}{x_4-\la}dx_4
 \end{align*}
  \item[$\bf(B_4)$]
$\boxed{(\la_1-\la_3)(\la_1-\la_4)(f_{13}f_{24}-f_{23}f_{14})+(\la_3-\la_4)(f_{23}f_{24}-f_{22}f_{34})=0}$
\begin{align*}
X_1(\la)=\frac{f_{24}(\la_1-\la_4)}{f_{22}}(\la_1-\la)\partial_{ x_1}-\frac{f_{24}(\la_4-\la)+f_{14}(\la_1-\la_4)(\la_1-\la)}{f_{22}}\partial_{ x_2}+
(\la_4-\la)\d_{ x_4}\\
X_2(\la)=\frac{f_{23}(\la_1-\la_3)}{f_{22}}(\la_1-\la)\partial_{ x_1}-\frac{f_{23}(\la_3-\la)+f_{13}(\la_1-\la_3)(\la_1-\la)}{f_{22}}\partial_{ x_2}
+(\la_3-\la)\d_{ x_3}\\
[X_1(\la),X_2(\la)]=0
 \end{align*}
 \begin{align*}
-f_2f_{23}+(\la_1-\la_3)(f_1 f_{23}- f_2f_{13})+f_{3}f_{22}=0\\
-f_2f_{24}+(\la_1-\la_4)(f_1 f_{24}- f_2f_{14})+f_{4}f_{22}=0\\
(\la_3-\la_4)(\la_1f_1-f_2)f_{34}+(\la_1-\la_3)\la_1f_4f_{13}-(\la_1-\la_4)\la_1f_3f_{14}
+\la_3f_4f_{23}-\la_4f_3f_{24}=0\\
(\la_1-\la_3)f_4f_{23}-(\la_1-\la_4)f_3f_{24}+(\la_3-\la_4)f_2f_{34}=0
 \end{align*}
  \begin{align*}
\al_\la=\frac{(\la_1-\la)f_1-f_2}{(\la_1-\la)^2}dx_1+
\frac{f_2}{\la_1-\la}dx_2+\frac{f_3}{\la_3-\la}dx_3+
\frac{f_4}{\la_4-\la}dx_4
 \end{align*}
  \begin{align*}
\om:=(\la_1-\la)^2d\al_\la\wedge dx_3\wedge dx_4=f_{22}dx_1\wedge dx_2\wedge dx_3\wedge dx_4
 \end{align*}
$$
g:=\frac{\la_3-\la_4}{J}\left[\begin{array}{cccc}
2f_{22}f_{14}f_{13} & f_{22}(f_{14}f_{23}+f_{13}f_{24}) & f_{22}f_{34}f_{13}+f_{14}f_{23}^2 & f_{22}f_{34}f_{14}+f_{13}f_{24}^2\\
* & 2 f_{22}f_{24}f_{23} & f_{23}(f_{23}f_{24}+f_{22}f_{34}) & f_{24}(f_{23}f_{24}+f_{22}f_{34}) \\
* & * & 2f_{34}f_{23}^2 & 2f_{24}f_{34}f_{23}\\
* & * & * & 2f_{34}f_{24}^2
\end{array}\right],
$$
where $J=f_{14}f_{23}-f_{13}f_{24}$.
\item[$\bf(C_0)$]
\framebox{\parbox{16.4cm}{$\color{black}((x_3-x_4)(x_2f_{22}-f_{12}+2f_2)+f_{22})f_{14}+
(2f_2-f_{12}-(x_3-x_4)(f_{23}+x_2f_{12}-f_{11}+f_1))f_{24}+(x_3-x_4)f_{22}f_{34}=0$}}
\begin{align*}
X_1(\la)=\frac{(x_3-\la)(x_3-x_4)f_{24}}{f_{22}}\partial_{ x_1}
-\frac{(x_3-\la)(x_3-x_4)f_{14}+(x_4-\la)f_{24}}{f_{22}}\partial_{ x_2}+
(x_4-\la)\d_{ x_4}\\
X_2(\la)=\frac{(x_3-\la)(x_2f_{22}-f_{12}+2f_2){\color{black}+f_{22}}}{f_{22}}\partial_{ x_1}+
\frac{(x_3-\la)(-x_2f_{12}+f_{11}-f_{23}-f_1)-f_{12}+2f_2}{f_{22}}\partial_{ x_2}\\
 +(x_3-\la)\partial_{ x_3}\\
 [X_1(\la),X_2(\la)]=0
 \end{align*}
 \begin{align*}
(-x_2f_2-f_1)f_{12}+(x_2f_1+f_3)f_{22}-f_2(f_{23}-f_1-f_{11})=0\\
(-f_2+(x_3-x_4)f_1)f_{24}-f_2(x_3-x_4)f_{14}+f_4f_{22}=0\\
{\color{black}(f_2-x_3f_1+(x_3-x_4)(x_2f_2+x_3f_3))f_{14}}-(x_3-x_4)(f_1x_3-f_2)f_{34}-(x_2x_3+1)f_4f_{12}\\
+x_4(x_2f_1+f_3)f_{24}+f_4(x_3(f_{11}-f_{23}-f_1)+2f_2)=0\\
(x_2x_3+1)f_2f_{24}+(x_3-x_4)(f_3x_3-f_1)f_{24}-(x_2x_3+1)f_4f_{22}\\
-x_3(x_3-x_4)f_2f_{34}+x_3f_4f_{12}-f_2(x_4f_{14}+2x_3f_4)=0
 \end{align*}
  \begin{align*}
\al_\la=\frac{(x_3-\la)f_1-f_2}{(x_3-\la)^2}dx_1
+\frac{f_2}{x_3-\la}dx_2
+\frac{(x_2(x_3-\la)+1)f_2+(x_3-\la)^2f_3-(x_3-\la)f_1}{(x_{\color{black}3}-\la)^3}dx_3\\
+\frac{f_4}{x_4-\la}dx_4
 \end{align*}
 \item[$\bf(C_3)$]
$\boxed{(\la_3-\la_4)(f_{14}f_{33}-f_{23}f_{24}+(-f_{13}+f_{22})f_{34})-f_{23}f_{34} +f_{24}f_{33}=0}$
\begin{align*}
X_1(\la)=\frac{(\la_3-\la)(\la_3-\la_4)f_{34}}{f_{33}}\partial_{ x_2}
-\frac{(\la_3-\la)(\la_3-\la_4)f_{24}+(\la_4-\la)f_{34}}{f_{33}}\d_{ x_3}
+(\la_4-\la)\d_{x_4}\\
X_2(\la)=(\la_3-\la)\d_{x_1}+\left(1-\frac{(\la_3-\la)f_{23}}{f_{33}}
\right)\partial_{ x_2}+
\frac{(\la_3-\la)(f_{22}-f_{13})-f_{23}}{f_{33}}\partial_{ x_3}\\
[X_1(\la),X_2(\la)]=0
\end{align*}
 \begin{align*}
f_3(f_{22}-f_{13})-f_2f_{23}+f_1f_{33}=0\\
(\la_3-\la_4)(f_2f_{34}-f_3f_{24})+f_4f_{33}-f_3f_{34}=0\\
((-f_2+f_1(\la_3-\la_4))\la_3+f_3)f_{24}-(f_2\la_3-f_3)f_{14}(\la_3-\la_4)
+f_4((f_{22}-f_{13})\la_3-f_{23})+f_1f_{34}\la_4=0\\
((\la_3-\la_4)(\la_3f_1-f_2)+f_3)f_{34}-\la_3(\la_3-\la_4)f_3f_{14}+\la_3f_4f_{23}-\la_4f_3f_{24}-f_4f_{33}=0
 \end{align*}
  \begin{align*}
\al_\la=-\frac{-(\la_3-\la)^2f_1+(\la_3-\la)f_2-f_3}{(\la_3-\la)^3}dx_1
+\frac{(\la_3-\la)f_2-f_3}{(\la_3-\la)^2}dx_2
+\frac{f_3}{\la_3-\la}dx_3+
\frac{f_4}{\la_4-\la}dx_4
 \end{align*}
  \begin{align*}
\om:=(\la_3-\la)^2d\al_\la\wedge dx_1\wedge dx_4=f_{33}dx_1\wedge dx_2\wedge dx_3\wedge dx_4
 \end{align*}
 $$
 g:=-\frac{\la_3-\la_4}{2J}\left[\begin{array}{cccc}
-2(f_{33}(f_{13}-f_{22})+f_{23}^2)f_{14} & b & c & d\\
* & 2 f_{33}f_{23}f_{24} & f_{33}(f_{34}f_{23}+f_{24}f_{33}) & f_{34}(f_{34}f_{23}+f_{24}f_{33}) \\
* & * & 2f_{33}^2f_{34} & 2f_{33}f_{34}^2\\
* & * & * & 2f_{34}^3
\end{array}\right],
 $$
 where $J=f_{34}f_{23}-f_{24}f_{33}$,
 \begin{align*}
b=-(f_{33}(f_{13}-f_{22})+f_{23}^2)f_{24}+f_{14}f_{23}f_{33},\\
c=(-f_{14}f_{33}+(-f_{13}+f_{22})f_{34})f_{33}+f_{23}^2f_{34},\\
d=(-f_{14}f_{33}-2f_{23}f_{24}+(-f_{13}+f_{22})f_{34})f_{34}+f_{24}^2f_{33}.
 \end{align*}
 \item[$\bf(D_0)$]
$\boxed{(x_1-x_3)^2(f_{13}f_{24}-f_{23}f_{14})+f_{24}^2-(f_{22}-f_{2})(f_{44}-f_4)=0}$
\begin{align*}
X_1(\la)=\frac{f_{24}(x_1-x_3)}{f_{22}-f_2}(x_1-\la)\partial_{ x_1}-\frac{f_{14}(x_1-x_3)(x_1-\la)+f_{24}(x_3-\la)}{f_{22}-f_2}\partial_{ x_2}+
(x_3-\la)\d_{ x_4}\\
X_2(\la)=\frac{f_{23}(x_1-x_3)-f_{24}}{f_{22}-f_2}(x_1-\la)\partial_{ x_1}-\frac{(x_3-\la)f_{23}+f_{13}(x_1-x_3)(x_1-\la)-f_{14}(x_1-\la)+f_{24}}{f_{22}-f_2}\partial_{ x_2}\\
+(x_3-\la)\d_{ x_3}+\d_{ x_4}\\
[X_1(\la),X_2(\la)]={\color{black}-}X_1(\la)
 \end{align*}
  \begin{align*}
-x_3(-f_2f_{23}-(x_1-x_3)(f_1 f_{23}- f_2f_{13}))-x_1f_2f_{14}-(x_3 f_3-f_4)f_{22}\\
+(x_1 f_1-f_2)f_{24}+f_2(x_3f_3-f_4)=0\\
-x_1(-f_4f_{14}-(x_1-x_3)(f_3 f_{14}- f_4f_{13}))-x_3f_4f_{23}-(x_3 f_3-f_4)f_{24}\\
+(x_1 f_1-f_2)(f_4-f_{44})=0\\
(-f_2+(x_1-x_3)f_1) f_{24}-(x_1-x_3) f_2f_{14}+f_4(f_{22}-f_2)=0\\
(-f_4-(x_1-x_3)f_3) f_{24}+(x_1-x_3) f_4f_{23}+f_2(f_{44}-f_4)=0
 \end{align*}
  \begin{align*}
\al_\la=\frac{(x_1-\la)f_1-f_2}{(x_1-\la)^2}dx_1+
\frac{f_2}{x_1-\la}dx_2+\frac{(x_3-\la)f_3-f_4}{(x_3-\la)^2}dx_3+
\frac{f_4}{x_3-\la}dx_4
 \end{align*}
  \item[$\bf(D_3)$]
$\boxed{(\la_1-\la_3)^2(f_{13}f_{24}-f_{23}f_{14})+f_{24}^2-f_{22}f_{44}=0}$
\begin{align*}
X_1(\la)=\frac{f_{24}(\la_1-\la_3)}{f_{22}}(\la_1-\la)\partial_{ x_1}-\frac{f_{14}(\la_1-\la_3)(\la_1-\la)+f_{24}(\la_3-\la)}{f_{22}}\partial_{ x_2}+
(\la_3-\la)\d_{ x_4}\\
X_2(\la)=\frac{f_{23}(\la_1-\la_3)-f_{24}}{f_{22}}(\la_1-\la)\partial_{ x_1}-\frac{(\la_3-\la)f_{23}+f_{13}(\la_1-\la_3)(\la_1-\la)-f_{14}(\la_1-\la)+f_{24}}{f_{22}}\partial_{ x_2}\\
+(\la_3-\la)\d_{ x_3}+\d_{ x_4}\\
[X_1(\la),X_2(\la)]=0
 \end{align*}
 \begin{align*}
\la_3(-f_2f_{23}+(\la_1-\la_3)(f_1 f_{23}- f_2f_{13}))+\la_1f_2f_{14}+(\la_3 f_3-f_4)f_{22}\\
-(\la_1 f_1-f_2)f_{24}=0\\
\la_1(-f_4f_{14}-(\la_1-\la_3)(f_3 f_{14}-f_4f_{13}))+\la_3f_4f_{23}-(\la_3 f_3-f_4)f_{24}\\
+(\la_1 f_1-f_2)f_{44}=0\\
(-f_2+(\la_1-\la_3)f_1) f_{24}-(\la_1-\la_3) f_2f_{14}+f_4f_{22}=0\\
(-f_4-(\la_1-\la_3)f_3) f_{24}+(\la_1-\la_3) f_4f_{23}+f_2f_{44}=0
 \end{align*}
  \begin{align*}
\al_\la=\frac{(\la_1-\la)f_1-f_2}{(\la_1-\la)^2}dx_1+
\frac{f_2}{\la_1-\la}dx_2+\frac{(\la_3-\la)f_3-f_4}{(\la_3-\la)^2}dx_3+
\frac{f_4}{\la_3-\la}dx_4
 \end{align*}
 \begin{align*}
\om:=(\la_1-\la)^2d\al_\la\wedge dx_3\wedge dx_4=f_{22}dx_1\wedge dx_2\wedge dx_3\wedge dx_4
 \end{align*}
$$
g:=\frac1{J}\left[\begin{array}{cccc}
2f_{22}f_{14}^2 & 2f_{22}f_{24}f_{14} & -f_{13}(f_{24}^2-f_{22}f_{44})+2f_{24}f_{14}f_{23} & f_{14}(f_{24}^2+f_{22}f_{44})\\
* & 2 f_{22}f_{24}^2 & f_{23}(f_{24}^2+f_{22}f_{44}) & f_{24}(f_{24}^2+f_{22}f_{44}) \\
* & * & 2f_{44}f_{23}^2 & 2f_{44}f_{23}f_{24}\\
* & * & * & 2f_{44}f_{24}^2
\end{array}\right],
$$
where $J=f_{14}f_{23}-f_{13}f_{24}$.
  \item[$\bf(E_0)$]
\framebox{\parbox{16.4cm}{$x_2(-f_{13}f_{23}+f_{33}f_{12})+2x_3(f_{22}f_{33} -f_{23}^2)+(f_{13}-3f_3)(f_{13}-f_{22})+f_{23}(-f_{34}+f_{12}-
2f_2)+f_{33}(f_{24}-f_{11}+f_1)=0$}}
{\color{black}
 \begin{align*}
X_1(\la)=\partial_{x_1}+\left(\frac{(x_4-\la)(x_2f_{23}+2x_3f_{33}-f_{13}+3f_3)}{f_{33}}-x_2\right)
\partial_{x_2}\\
+\left(\frac{(\la-x_4)(f_{24}+x_2 f_{12}-f_{11}+2x_3 f_{22}+f_1)}{f_{23}}\right.\\
+\left.\frac{(x_2f_{23}-f_{13}+3f_3)((\la-x_4)(f_{22}-f_{13})+f_{23})}{f_{33}f_{23}}\right)
\partial_{x_3}+(x_4-\la)\partial_{x_4}\\
X_2(\la)=(x_4-\la)\partial_{x_1}+\frac{(\la-x_4)f_{23}+f_{33}}{f_{33}}\partial_{x_2}
-\frac{(\la-x_4)(f_{22}-f_{13})+f_{23}}{f_{33}}\partial_{x_3}\\
[X_1(\la),X_2(\la)]=X_2(\la)
 \end{align*}
 }
 \begin{align*}
f_{33}f_1+f_{22}f_3-f_{13}f_3-f_{23}f_2=0\\
((2x_3x_4-x_2)f_1-2x_3f_2-f_4)f_{23}+(-x_4f_1-x_2(x_4f_2-f_3))f_{12}+(-2x_3x_4f_2-x_4f_4+f_1)f_{13}\\
+(x_2x_4f_1+2x_3f_3+x_4f_4)f_{22}+(x_4f_2-f_3)f_{11}+(-x_4f_2+f_3)f_{24}+f_1(x_4f_{34}+x_4f_2-2f_3)=0\\
((-2x_3x_4+x_2)f_3-x_4f_1+f_2)f_{13}+((2x_3x_4-x_2)f_1-2x_3f_2-f_4)f_{33}+(x_2x_4f_1+2x_3f_3+x_4f_4)f_{23}\\
+((-x_2x_4-1)f_{12}+x_4(f_{11}-f_{24}+2f_1)+f_{34}-f_2)f_3=0\\
(x_2x_4f_2-2x_3x_4f_3+f_2)f_{23}+(2x_3x_4f_2+x_4f_4-f_1)f_{33}+(-x_4f_2+f_3)f_{13}\\
+((-x_2x_4-1)f_{22}+x_4(f_{12}-f_{34}+f_2))f_3=0
 \end{align*}
  \begin{align*}
\al_\la=((x_4-\la)^2f_1-(x_4-\la)f_2+f_3)\frac{dx_1}{(x_{\color{black}4}-\la)^3}
+((x_4-\la)f_2-f_3)\frac{dx_2}{(x_{\color{black}4}-\la)^2}+f_3\frac{dx_3}{x_{\color{black}4}-\la}\\
+((x_4-\la)^3f_4+(x_4-\la)^2(2x_3f_3+x_2f_2-f_1)+(x_4-\la)(-x_2f_3+f_2)-f_3)\frac{dx_4}{(x_{\color{black}4}-\la)^4}
 \end{align*}
 \item[$\bf(E_1)$]
$\boxed{(f_{14}-f_{23})f_{34}+(f_{22}-f_{13})f_{44}+(f_{33}-f_{24})f_{24}=0}$
\begin{align*}
X_1(\la)=(\la_4-\la)\partial_{ x_1}+\partial_{ x_2}-\frac{f_{24}(\la_4-\la)}{f_{44}}\partial_{ x_3}+
\frac{(\la_4-\la)(f_{23}-f_{14})-f_{24}}{f_{44}}\d_{ x_4}\\
X_2(\la)=(\la_4-\la)\partial_{ x_2}+\frac{f_{44}-f_{34}(\la_4-\la)}{f_{44}}\partial_{ x_3}+
\frac{(\la_4-\la)(f_{33}-f_{24})-f_{34}}{f_{44}}\d_{ x_4}\\
[X_1(\la),X_2(\la)]=0
 \end{align*}
 \begin{align*}
(f_3-\la_4 f_2)f_{24}+\la_4f_4(f_{22}-f_{13})+ f_1{\color{black}(\la_4f_{34}-f_{44})}+f_4(f_{14}-f_{23})=0\\
(f_4-\la_4 f_3)f_{24}+\la_4f_4(f_{23}-f_{14})+ (\la_4 f_1-f_2)f_{44}-f_4f_{33}+f_3f_{34}=0\\
{\color{black}(f_{22}-f_{13})(\la_4 f_{3}-f_4)}+(-\la_4 f_1+f_2)f_{24}+\la_4 f_2(f_{14}-f_{23})+f_1(\la_4 f_{33}-f_{34})=0\\
f_2f_{44}-f_3f_{34}+f_4(f_{33}-f_{24})=0
 \end{align*}
 \begin{align*}
\al_\la=\frac{(\la_4-\la)^3f_1-(\la_4-\la)^2f_2+(\la_4-\la)f_3-f_4}{(\la_4-\la)^4}dx_1+
\\
\frac{(\la_4-\la)^2f_2-(\la_4-\la)f_3+f_4}{(\la_4-\la)^3}dx_2
+\frac{(\la_4-\la)f_3-f_4}{(\la_4-\la)^2}dx_3+
\frac{f_4}{\la_4-\la}dx_4
 \end{align*}
 $$
 \om=(\la_4-\la)^2d\al_\la\wedge dx_1\wedge dx_2=f_{44}dx_1\wedge dx_2\wedge dx_3\wedge dx_4
 $$
$$
g:=\frac1{J}\left[\begin{array}{cccc}
a & b & c & d\\
* & 2f_{24}(f_{34}^2+f_{44}(f_{24}-f_{33})) & f_{34}(f_{34}^2-f_{44}(f_{33}-2f_{24})) & f_{44}(f_{34}^2-f_{44}(f_{33}-2f_{24})) \\
* & * & 2f_{44}f_{34}^2 & 2f_{44}^2f_{34}\\
* & * & * & 2f_{44}^3
\end{array}\right],
$$
where $J=f_{33}f_{44}-f_{34}^2$,
\begin{align*} a=2[((f_{13}-f_{22})f_{24}+f_{13}f_{33}+f_{14}^2-2f_{23}f_{14}-f_{22}f_{33}+f_{23}^2)f_{44}+f_{24}^3\\
+f_{34}f_{24}(f_{14}-f_{23})-f_{34}^2(f_{13}-f_{22})], \\ b=(2f_{24}-f_{13})(f_{14}-f_{23})f_{44}+f_{34}((f_{14}-f_{23})f_{34}+2f_{24}^2),\\
c=(2f_{24}-f_{33})f_{34}^2+2f_{34}f_{44}(f_{14}-f_{23})+f_{44}f_{33}^2,\\
d=2(f_{14}-f_{23})f_{44}^2+f_{34}f_{44}(2f_{24}+f_{33})+f_{34}^3.
 \end{align*}
 \item[$\bf(F_0)$]
$\boxed{(x_1-x_3)(x_1-x_4)(f_{14}f_{23}-f_{13}f_{24})-(x_3-x_4)f_2f_{34}=0}$
\begin{align*}
X_1(\la)=-\frac{f_{24}(x_1-x_4)}{f_2}(x_1-\la)\partial_{ x_1}+\frac{f_{14}(x_1-x_4)}{f_2}(x_1-\la)\partial_{ x_2}+
(x_4-\la)\d_{ x_4}\\
X_2(\la)=-\frac{f_{23}(x_1-x_3)}{f_2}(x_1-\la)\partial_{ x_1}+\frac{f_{13}(x_1-x_3)}{f_2}(x_1-\la)\partial_{ x_2}+
(x_3-\la)\d_{ x_3}\\
[X_1(\la),X_2(\la)]=0
 \end{align*}
 \begin{align*}
(x_1-x_3)(f_1f_{23}-f_2f_{13})-f_2f_3=0\\
(x_1-x_4)(f_1f_{24}-f_2f_{14})-f_2f_4=0\\
(x_1-x_3)f_4f_{13}-(x_1-x_4)f_3f_{14}+(x_3-x_4)f_1f_{34}=0\\
(x_1-x_3)f_4f_{23}-(x_1-x_4)f_3f_{24}+(x_3-x_4)f_2f_{34}=0
 \end{align*}
  \begin{align*}
\al_\la=\frac{f_1}{x_1-\la}d{x_1}+
\frac{f_2}{x_1-\la}d{x_2}+\frac{f_3}{x_3-\la}d{x_3}+
\frac{f_4}{x_4-\la}d{x_4}
 \end{align*}
  \item[$\bf(G_0)$]
$\boxed{(x_1-x_3)^2(f_{14}f_{23}-f_{13}f_{24})+f_2f_{4}=0}$
\begin{align*}
X_1(\la)=-\frac{f_{23}(x_1-x_3)}{f_2}(x_1-\la)\partial_{ x_1}+\frac{f_{13}(x_1-x_3)}{f_2}(x_1-\la)\partial_{ x_2}+
(x_3-\la)\d_{ x_3}\\
X_2(\la)=-\frac{f_{24}(x_1-x_3)}{f_2}(x_1-\la)\partial_{ x_1}+\frac{f_{14}(x_1-x_3)}{f_2}(x_1-\la)\partial_{ x_2}+
(x_3-\la)\d_{ x_4}\\
[X_1(\la),X_2(\la)]=X_2(\la)
 \end{align*}
 \begin{align*}
(x_1-x_3)(f_1f_{23}-f_2f_{13})-f_2f_3=0\\
(x_1-x_3)(f_1f_{24}-f_2f_{14})-f_2f_4=0\\
(x_1-x_3)(f_4f_{13}-f_3f_{14})-f_1f_{4}=0\\
(x_1-x_3)(f_4f_{23}-f_3f_{24})-f_2f_{4}=0
 \end{align*}
  \begin{align*}
\al_\la=\frac{f_1}{x_1-\la}d{x_1}+
\frac{f_2}{x_1-\la}d{x_2}+\frac{f_3}{x_3-\la}d{x_3}+
\frac{f_4}{x_3-\la}d{x_4}
 \end{align*}
  \item[$\bf(H_0)$]
$\boxed{(x_1-x_3)^2(f_{14}f_{23}-f_{13}f_{24})-f_{4}(f_{22}-f_2)=0}$
\begin{align*}
X_1(\la)=-\frac{f_{23}(x_1-x_3)}{f_2-f_{22}}(x_1-\la)\partial_{ x_1}+\frac{f_{13}(x_1-x_3)(x_1-\la)+f_{23}(x_3-\la)}{f_2-f_{22}}\partial_{ x_2}+
(x_3-\la)\d_{ x_3}\\
X_2(\la)=-\frac{f_{24}(x_1-x_3)}{f_2-f_{22}}(x_1-\la)\partial_{ x_1}+\frac{f_{14}(x_1-x_3)(x_1-\la)+f_{24}(x_3-\la)}{f_2-f_{22}}\partial_{ x_2}+
(x_3-\la)\d_{ x_4}\\
[X_1(\la),X_2(\la)]=X_2(\la)
 \end{align*}
 \begin{align*}
(x_1-x_3)(f_1f_{23}-f_2f_{13})+f_3(f_{22}-f_2)-f_2f_{23}=0\\
(x_1-x_3)(f_1f_{24}-f_2f_{14})+f_4(f_{22}-f_2)-f_2f_{24}=0\\
x_1(x_1-x_3)(f_4f_{13}-f_3f_{14})+x_3(f_4f_{23}-f_3f_{24})-(x_1f_1-f_2)f_4=0\\
(x_1-x_3)(f_4f_{23}-f_3f_{24})-f_2f_{4}=0
 \end{align*}
  \begin{align*}
\al_\la=\frac{(x_1-\la)f_1-f_2}{(x_1-\la)^2}d{x_1}+
\frac{f_2}{x_1-\la}d{x_2}+\frac{f_3}{x_3-\la}d{x_3}+
\frac{f_4}{x_3-\la}d{x_4}
 \end{align*}
 \end{enumerate}

\section{Universal hierarchy and the corresponding heavenly equation}
\label{hyperCR}

\abz\label{theoHCR}
\begin{theo} \cite[Ch. I, Prop. 5]{rigal} Let $\al_\la=\al_0+\la\al_1+\cdots+\la^n\al_n$ be a $\la$-depending one-form defining a Veronese web on a $(n+1)$-dimensional manifold $M$, i.e. the one-forms $\al_i$ are linearly independent and the Frobenius integrability condition
\begin{equation}\label{frob}\equ
\al_\la\wedge d\al_\la=0
\end{equation}
holds. Then in a neighbourhood of any point in $M$ there exists a local coordinate system $(x_0,\ldots,x_n)$ and a function $g$ such that $\al_\la$ is proportional to the form
\begin{equation}\label{frob1}\equ
dx_0+\la(dx_1-g_{x_n}dx_0)+\la^2(dx_2-g_{x_n}dx_1-g_{x_{n-1}}dx_0)+
\cdots+\la^n(dx_n-g_{x_n}dx_{n-1}+\cdots-g_{x_1}dx_0).
\end{equation}
\end{theo}

\abz\label{coroHCR}
\begin{coro} Frobenius integrability condition (\ref{frob}) is equivalent to a system of second order PDEs on the function $g$. For $n=2$  this system is given by a single equation
\begin{equation}\label{eHCR}\equ
g_2g_{12}-g_1g_{22}+g_{11}-g_{02}=0
\end{equation}
and {\color{black}for $n=3$ it is of the form}
\begin{align}\label{hcr}\equ
\begin{split}
g_2g_{13}-g_1g_{23}+g_{11}-g_{02}&=0\\
g_3g_{13}-g_1g_{33}+g_{12}-g_{03}&=0\\
g_3g_{23}-g_2g_{33}+g_{22}-g_{13}&=0\\
g_3g_{02}-g_3g_{11}+g_2g_{12}-g_2g_{03}+g_1g_{13}-g_1g_{22}&=0.
\end{split}
\end{align}

\end{coro}
\noindent{\sc Proof} Condition (\ref{frob}) is equivalent to
\begin{equation}\label{frob2}\equ
\sum_{i=0}^k\al_i\wedge d\al_{k-i}=0, \qquad k=0,\ldots,2n,
\end{equation}
where we put $\al_i=0$ for $i>n$.

Now the proof is a direct check of the restrictions on the function $g$ imposed by (\ref{frob2}) with $k>n+1$. \qed

\abz\label{remHCR}
\begin{rema}\rm
Equation (\ref{eHCR}) up to a shift of indices coincides with the hyper-CR equation \cite{dunajski}, \cite{dunajskiKr} decribing Veronese webs in 3D, {\color{black}which is not covered by the formalism of \cite{pHirota}. In \cite{KrynskiDef} this equation is obtained from a generic homogeneous Hirota PDE by a series of coordinate changes and limiting procedures. }

First three equations of system (\ref{hcr}) up to a permutation of indices coincide with the first three equations of the universal hierarchy of Martinez–Alonso and Shabat \cite{martinezAlShabat}, \cite{ferapontovKruglikov}.  {\color{black}The higher dimensional generalizations of equations (\ref{hcr}) as well as the 1-form given by (\ref{frob1}) were also indicated in \cite{dunajski}, \cite{ferapontovKruglikov}.}
\end{rema}

\abz\label{theo2HCR}
\begin{theo} The first three equations of system (\ref{hcr}) imply the fourth one {\color{black}(cf. footnote \ref{ft43} on page \pageref{ft43})}. If a function $g$ satisfies system (\ref{hcr}), then it satisfies the following "heavenly" PDE:
\begin{equation}\label{uni}\equ
-g_{02}g_{33}+g_{03}g_{23}+g_{11}g_{33}-g_{12}g_{23}-g_{13}^2+g_{13}g_{22}=0,
\end{equation}
{\color{black} which is equivalent to condition $d\al_\la\wedge d\al_\la=0$}.
This equation up to a shift of indices coincides with equation $\mathbf{E_1}$.
\end{theo}
\noindent{\sc Proof} {\color{black}Check that  $-g_3(I)+g_2(II)-g_1(III)$,  where (I)--(III) denote the first three equations of system (\ref{hcr}), gives the fourth equation. To prove the second statement one notices that (\ref{frob}) obviously implies $d\al_\la\wedge d\al_\la=0$, or alternatively checks that  $g_{33}(I)-g_{23}(II)+g_{13}(III)$ gives (\ref{uni})}. \qed

\section{Some examples and discussion of existence of  {\color{black} Mason--Newman vector fields} for $\mathbf{A_0}$  and $\mathbf{G_0}$ cases}
\label{sexa}

\abz\label{exaI}
\begin{exa}\rm
Consider the function $f=2(x_2-x_3)(x_2-x_4)/(x_1-x_2)$. It is a solution to the Hirota system and the heavenly equation corresponding to $\mathbf{A_0}$ case. The Mason--Newman vector fields
\begin{align*}
 X_1(\lambda)&=\frac{x_{13}x_{24}}{q(x)}(x_1-\lambda)\partial_{x_1}+\frac{x_{14}x_{23}}{q(x)}(x_2-\lambda)\partial_{x_2}
 +(x_4-\lambda)\partial_{x_4}\\
 X_2(\lambda)&=\frac{x_{14}x_{23}}{q(x)}(x_1-\lambda)\partial_{x_1}+\frac{x_{13}x_{24}}{q(x)}(x_2-\lambda)\partial_{x_2}
 +(x_3-\lambda)\partial_{x_3},
\end{align*}
where $q(x):=x_{13}x_{24}+x_{14}x_{23}$ and
$$
x_{ij}:=x_i-x_j,
$$
 are divergence free with respect to the volume form
$$
\omega = \frac{q(x)x_{12}^2}{(x_{13}x_{14}x_{23}x_{34})^2}dx_1\wedge dx_2\wedge dx_3\wedge dx_4.
$$
The corresponding self-dual vacuum Einstein metric is {\em flat} and is given by the following matrix
$$
\frac{x_{12}}{J(x)}\begin{bmatrix}
 q(x)x_{23}x_{24} & -\frac{1}{2}q(x)^2 & \frac{1}{2}\big(q(x)+x_{14}x_{23}\big)x_{12}x_{24} & \frac{1}{2}\big(q(x)+x_{13}x_{24}\big)x_{12}x_{23} \\
  \ast & q(x)x_{13}x_{14} & -\frac{1}{2}\big(q(x)+x_{13}x_{24}\big)x_{12}x_{14} & -\frac{1}{2}\big(q(x)+x_{14}x_{23}\big)x_{12}x_{13}\\
 \ast & \ast &  x_{12}^2x_{14}x_{24} & \frac{1}{2}q(x)x_{12}^2 \\
 \ast & \ast & \ast & x_{12}^2x_{13}x_{23}
\end{bmatrix},
$$
where $J(x)=(x_{13}x_{14}x_{23}x_{34})^2$.

\end{exa}

\abz\label{remaflat}
\begin{rema}\rm
It is interesting to observe that the Veronese web corresponding to the above solution of the Hirota system is non-flat. A Veronese web $\{\F_\la\}$ is said to be {\em flat} if in a vicinity of any point there exists a local diffeomorphism to an open set in a vector space bringing $\F_\la$ to the foliation of parallel hyperplanes. If a Veronese web is defined by a 1-form $\al_\la$ satisfying (\ref{frob}) and consequently (\ref{frob2}), it is flat if and only if \cite[Ch. I, Prop. 6]{rigal}
\begin{equation*}
d\al_1 \wedge \al_1=0.
\end{equation*}
One easily checks that this condition is not satisfied for the Veronese web {\color{black} corresponding to the  above solution of the Hirota system}, which gives an example of including of a flat divergence free Kronecker web of rank 2 into a non-flat Veronese web of rank 3, where by the flatness of the former we mean the  flatness of the corresponding Einstein metric.
\end{rema}

\abz\label{remas}
\begin{rema}\rm
In frames of Example \ref{exaI}, consider the 2-form $\be_\la=i(X_1(\la)\wedge X_2(\la))\om$ annihilating $\langle X_1(\la),X_2(\la)\rangle$ and the $2$-form
$\be=i(X_1(\infty)\wedge X_2(\infty))\om=\be_\infty$ which was used in Theorem \ref{th4} {\color{black}and Corollary
\ref{corrop}} for establishing of  the existence of the potential. {\color{black} The 2-form $\be$ } is $d$-closed, but it is not $d_N$-closed with respect to the Nijenhuis operator $N:\partial_{x_i}\mapsto x_i \partial_{x_i}$ extending the partial Nijenhuis operator $X_i(\infty)\mapsto X_i(0)$, $i=1,2$,  since the vector fields $X_1(\infty),X_2(\infty)$ do not annihilate the function $\Tr(N)$. Note however, that there exists another $\la$-depending $2$-form $\tilde{\be}_\la=dd_{(N-\la I)^{-1}}f$ annihilating $\langle X_1(\la),X_2(\la)\rangle$ and its limit at infinity, $\tilde{\be}_\infty=dd_{N}f$ (cf. Theorem \ref{th561}), is obviously $d$- and $d_N$-closed (the dependence of $\tilde{\be}_\la$ on $\la$ is not polynomial). This justifies the existence of the potential in the example.
\end{rema}

The example above is very special and we claim that the Kronecker web defined by a general solution of heavenly equation $\mathbf{A_0}$ is not divergence free in the sense of Definition \ref{divfreedefi}. To argue this let $f$ be a solution of heavenly equation $\mathbf{A_0}$ and consider {\color{black}the corresponding closed $2$-form $(x_{ij}:=x_i-x_j)$
\begin{align*}
\tilde{\be}_\la=dd_{(N-\la I)^{-1}}f=\frac{f_{12}x_{12}}{(x_1-\la) (x_2-\la)}dx_1 \wedge dx_2+
\frac{f_{13}x_{13}}{(x_1-\la) (x_3-\la)}dx_1 \wedge dx_3\\
\frac{f_{14}x_{14}}{(x_1-\la) (x_4-\la)}dx_1 \wedge dx_4+
\frac{f_{23}x_{23}}{(x_2-\la) (x_3-\la)}dx_2 \wedge dx_3\\
+\frac{f_{24}x_{24}}{(x_2-\la) (x_4-\la)}dx_2 \wedge dx_4+\frac{f_{34}x_{34}}{(x_3-\la) (x_4-\la)}dx_3 \wedge dx_4
\end{align*}
 (here $N$ is given by the corresponding formula in Appendix). Assume there exist linearly independent vector fields $ X_1(\la),X_2(\la)$ annihilating $\be_\la$ and a volume form $\om$ with respect to which they  are divergence free. Then the $2$-form $\be_\la=i(X_1(\la)\wedge X_2(\la))\om$ is $d$-closed (see Theorem \ref{th1}), is a polynomial of second degree in $\la$,  and also annihilates $X_{1,2}(\la)$.} Hence $\be_\la$ has to be proportional to the $2$-form $(x_1-\la)(x_2-\la)(x_3-\la)(x_4-\la)\tilde\be_\la$, which is also a polynomial of the second degree in $\la$ and annnihilates {\color{black} $X_{1,2}(\la)$}, with {\color{black}a} coefficient of proportionality {\color{black}$U(x)$} not depending on $\la$:
$$
\be_\la=U(x)(x_1-\la)(x_2-\la)(x_3-\la)(x_4-\la)\tilde\be_\la.
$$
Condition of closedness of $\be_\la$ implies the following system of linear PDEs on the function $U$:
\begin{align*}
f_{23} x_{23}U_1-f_{13} x_{13}U_2+f_{12} x_{12}U_3=&0\\
f_{24}x_{24}U_1-f_{14} x_{14}U_2+f_{12} x_{12}U_4=&0\\
f_{23} x_{23} x_1U_1-f_{13}x_{13} x_2U_2+f_{12} x_{12} x_3U_3=&
- U(f_{12}x_{12}-f_{13}x_{13}+f_{23}x_{23})\\
f_{24} x_{24} x_1U_1-f_{14} x_{14} x_2U_2+f_{12} x_{12} x_4U_4=&
- U(f_{12}x_{12}-f_{14} x_{14}+f_{24}x_{24}).
\end{align*}
Assuming that $f_{13}f_{24}-f_{14}f_{23}\not=0$ one can solve this system with respect to $U_i$, getting
$U_i(x)=U(x)\cdot F_i(f_{kl}(x),x)$, $i=1,\ldots,4$, where the expressions $F_i$ depend on the second derivatives $f_{kl}(x)$ and explicitly on $x$. The necessary compatibility conditions implied by $U_{ij}(x)=U_{ji}(x)$ are given by
\begin{equation}\label{compaa}
\equ \partial_{x_j}F_i(f_{kl}(x),x)=\partial_{x_i}F_j(f_{kl}(x),x).
\end{equation}
One checks directly that these conditions are satisfied for the solution $f(x)$ from Example \ref{exaI}. On the other hand conditions (\ref{compaa}) are not satisfied for instance for the function
\begin{align*}
f(x)=\frac{x_{23}x_{24}x_{34}}{2x_{23}x_{14}-3x_{24}x_{13}+6x_{34}x_{12}}
\end{align*}
which is also a solution of heavenly equation $\mathbf{A_0}$; this finishes the proof of our claim.

\abz\label{exaII}
\begin{exa}\rm
It is worth to notice that {\color{black} still} there exists an example corresponding to $\mathbf{A_0}$ case leading to a non-flat self-dual vacuum Einstein metric.
The function $f=(x_1-x_4)^r (x_2-x_4)^r(x_3-x_4)^r$ satisfies Hirota system (\ref{HirotaA0}) for any $r\in \R$. The commuting vector fields $(x_{ij}:=x_i-x_j)$
\begin{align*}
 X_1(\la) = Y_1-\la X_1 =& -\frac{x_{14}x_{34}-r(x_{14}x_{24}+x_{14}x_{34}+x_{24}x_{34})}{rx_{12}x_{34}}(\la-x_1)\partial_{x_1}\\
 &+\frac{x_{24}x_{34}-r(x_{14}x_{24}+x_{14}x_{34}+x_{24}x_{34})}{rx_{12}x_{34}}(\la-x_2)\partial_{x_2}-(\la-x_4)\partial_{x_4}\\
 X_2(\la) =Y_2-\la X_2 =& -\frac{x_{23}x_{14}}{x_{12}x_{34}}(\la-x_1)\partial_{x_1}
 +\frac{x_{13}x_{24}}{x_{12}x_{34}}(\la-x_2)\partial_{x_2}-(\la-x_3)\partial_{x_3}
\end{align*}
are divergence free with respect to the  volume form
$$
    \omega = \frac{x_{12}}{x_{14}^2x_{24}^2x_{34}}dx_1\wedge dx_2\wedge dx_3\wedge dx_4,
$$
where we assumed $r\neq 0$.
The metric constructed {\color{black} by} Mason--Newman method is given by
$$
g=\begin{pmatrix}
\frac{r q_1}{x_{14}^2} & \frac{r(q_1+q_2)}{2x_{14}x_{24}} & \frac{r(q_1+q_3)}{2x_{14}x_{34}} & -\frac{q_1(q_2+q_3)}{2x_{14}} \\
\frac{r(q_1+q_2)}{2x_{14}x_{24}} & \frac{r q_2}{x_{24}^2} & \frac{r(q_2+q_3)}{2x_{24}x_{34}} & -\frac{q_2(q_1+q_3)}{2x_{24}} \\
\frac{r(q_1+q_3)}{2x_{14}x_{34}} & \frac{r(q_2+q_3)}{2x_{24}x_{34}} &
\frac{r q_3}{x_{34}^2} & -\frac{q_3(q_1+q_2)}{2x_{34}}  \\
-\frac{q_1(q_2+q_3)}{2x_{14}} & -\frac{q_2(q_1+q_3)}{2x_{24}} & -\frac{q_3(q_1+q_2)}{2x_{34}} & \frac{q_1q_2q_3}{r}
\end{pmatrix},
$$
where
$$
q_1=\frac{r-1}{x_{14}}+\frac{r}{x_{24}}+\frac{r}{x_{34}}, \qquad
q_2=\frac{r}{x_{14}}+\frac{r-1}{x_{24}}+\frac{r}{x_{34}}, \qquad
q_3=\frac{r}{x_{14}}+\frac{r}{x_{24}}+\frac{r-1}{x_{34}}.
$$
It is Ricci-flat and the only non-vanishing component of the Riemann tensor is given by
$$
R_{1313}=-\frac{(r-1)(2r-1)(3r-1)}{r^3x_{14}x_{24}x_{34}}.
$$
where the indices correspond to the {\color{black} coframe} $(X^1,X^2,Y^1,Y^2)$ dual to $(X_1,X_2,Y_1,Y_2)$.
\end{exa}

{\color{black} We conclude this section by a discussion showing that some solutions of  $\mathbf{G_0}$ heavenly equation can lead to self-dual vacuum Einstein metrics.}

It is not difficult to check that the  vector fields
\begin{align*}
X_1(\lambda)&=(x_1-\lambda)\partial_{x_1}+\frac{f_{14}}{f_{24}}(\lambda-x_1)\partial_{x_2}+\frac{(\lambda-x_3)f_2}{(x_1-x_3)f_{24}}\partial_{x_4}\\
 X_2(\lambda)&=-\frac{(\lambda-x_1)f_4}{(x_1-x_3)f_{24}}\partial_{x_2}+
 (x_3-\lambda)\partial_{x_3}+\frac{f_{23}}{f_{24}}(\lambda-x_3)\partial_{x_4}
\end{align*}
commute if the function $f$ satisfies {\color{black} $\mathbf{G_0}$ equation}
\begin{equation*}
(x_1-x_3)^2(f_{14}f_{23}-f_{13}f_{24})+f_2f_{4}=0. \label{1}
\end{equation*}
The necessary condition for the existence of a volume  form {\color{black} similar to (\ref{compaa})} implies that the function $f$  {\color{black} should} not depend on one of the variables.  Coordinate transformations lead to
\begin{align*}
X_1(\lambda) &= (y_1-\lambda)\partial_{y_1}+
\frac{(\lambda-y_2)y_3}{y_1-y_2}\left(\partial_{y_3}-F_3\partial_{y_4}\right)\\
X_1(\lambda) &= (y_2-\lambda)\partial_{y_2}-
\frac{(\lambda-y_1)y_3}{y_1-y_2}\partial_{y_3}+(\lambda-y_2)F_2 \partial_{y_4},
\end{align*}
where $F(y_1,y_2,y_3)$ satisfies the linear equation
\begin{equation}\label{leq}\equ
(y_1-y_2)^2\tilde{F}_{12}+y_3^2 \tilde{F}_{33} =0 , \qquad \tilde{F}=y_3F.
\end{equation}
{\color{black}The  corresponding volume form and the self-dual metric admitting the Killing vector $\partial_{y_4}$  are given by
\begin{align*}\label{metrlin}
\omega = -\frac{dy_1\wedge dy_2\wedge dy_3\wedge dy_4}{(y_1-y_2)y_3^2},\
g=\frac{2}{y_3^2F_3} \left(dy_4+F_2 dy_2)\odot(dy_4+F_2 dy_2+F_3 dy_3\right)-\frac{2F_3}{(y_1-y_2)^2}dy_1\odot dy_2.
\end{align*}}
Note that self-dual metrics which possess a Killing vector have been considered in \cite{FinPleb}, where also an equation similar to (\ref{leq})  on the key function appeared.

\section{Higher-dimensional generalizations}\label{6D}

Let $N:TM\to TM$ be a Nijenhuis operator, $f$ a function on $M$.  Retain assumptions from the beginning of Section \ref{secMain}. The distribution $D_\la:=\ker\be_\la$, where $\be_\la:=d\al_\la=d((N_\la)^{-1})^tdf$, is integrable as the kernel of a closed 2-form. One can achieve nontriviality of this distribution in many ways. In particular one of the possibilities is condition  (\ref{cond1}), or equivalently (\ref{cond2}), which means that the 2-form $\be_\la$ is of rank 2 and $D_\la$ is of rank $\dim M-2$. Starting from dimension 5 we get a system of PDEs on the function $f$. For instance, in dimension 5, taking $N$ given by $\d_{x_i}\mapsto x_i\d_{x_i}$ we get the system of equations
\begin{align*}
(x_1-x_2)(x_3-x_4)f_{12}f_{34}-(x_1-x_3)(x_2-x_4)f_{13}f_{24}
 +(x_1-x_4)(x_2-x_3)f_{14}f_{23}=0\\
 (x_1-x_2)(x_3-x_5)f_{12}f_{35}-(x_1-x_3)(x_2-x_5)f_{13}f_{25}
 +(x_1-x_5)(x_2-x_3)f_{15}f_{23}=0\\
 (x_1-x_2)(x_5-x_4)f_{12}f_{54}-(x_1-x_5)(x_2-x_4)f_{15}f_{24}
 +(x_1-x_4)(x_2-x_5)f_{14}f_{25}=0\\
 (x_1-x_5)(x_3-x_4)f_{15}f_{34}-(x_1-x_3)(x_5-x_4)f_{13}f_{54}
 +(x_1-x_4)(x_5-x_3)f_{14}f_{53}=0\\
 (x_5-x_2)(x_3-x_4)f_{52}f_{34}-(x_5-x_3)(x_2-x_4)f_{53}f_{24}
 +(x_5-x_4)(x_2-x_3)f_{54}f_{23}=0
\end{align*}
and the following parameter depending commuting vector fields spanning $D_\la$:
\begin{align*}
X_1(\la)=\frac{f_{25}(x_2-x_5)}{f_{12}(x_1-x_2)}(x_1-\la)\partial_{ x_1}-\frac{f_{15}(x_1-x_5)}{f_{12}(x_1-x_2)}(x_2-\la)\partial_{ x_2}+
(x_5-\la)\d_{ x_5}\\
X_2(\la)=\frac{f_{24}(x_2-x_4)}{f_{12}(x_1-x_2)}(x_1-\la)\partial_{ x_1}-\frac{f_{14}(x_1-x_4)}{f_{12}(x_1-x_2)}(x_2-\la)\partial_{ x_2}+
(x_4-\la)\d_{ x_4}\\
X_3(\la)=\frac{f_{23}(x_2-x_3)}{f_{12}(x_1-x_2)}(x_1-\la)\partial_{ x_1}-\frac{f_{13}(x_1-x_3)}{f_{12}(x_1-x_2)}(x_2-\la)\partial_{ x_2}+
(x_3-\la)\d_{ x_3}.
 \end{align*}
In turn, if, for instance, $N$ is the Jordan block with the constant eigenvalue $\la_1$,
$$
N:=\left[
           \begin{array}{ccccc}
             \la_1  & 0 & 0 & 0 & 0\\
             1 & \la_1 & 0 & 0 & 0\\
             0 & 1 & \la_1 & 0 & 0\\
              0 & 0 & 1 & \la_1 & 0\\
              0& 0 & 0 & 1 & \la_1
           \end{array}
         \right],
$$
we get
\begin{align*}
(f_{13}-f_{22})(f_{35}-f_{44})-(f_{14}-f_{23})(f_{25}-f_{34})+(f_{15}-f_{24})(f_{24}-f_{33})=0\\
(f_{13}-f_{22})f_{45}-(f_{14}-f_{23})f_{35}+(f_{24}-f_{33})f_{25}=0\\
(f_{13}-f_{22})f_{55}-(f_{15}-f_{24})f_{35}+(f_{25} -f_{34})f_{25}=0\\
(f_{14}-f_{23})f_{55}-(f_{15}-f_{24})f_{45}+(f_{35}-f_{44})f_{25}=0\\
(f_{24}-f_{33})f_{55}-(f_{25}-f_{34})f_{45}+(f_{35}-f_{44})f_{35}=0
\end{align*}
and
\begin{align*}
X_1(\la)&=(\la_1-\la)\partial_{x_1}+\partial_{x_2}-\frac{f_{25}}{f_{55}}(\la_1-\la)\partial_{x_4}+
\frac{(f_{24}-f_{15})(\la_1-\la)-f_{25}}{f_{55}}\partial_{x_5}\\
X_2(\la)&=(\la_1-\la)\partial_{x_2}+\partial_{x_3}-\frac{f_{35}}{f_{55}}(\la_1-\la)\partial_{x_4}+
\frac{(f_{25}-f_{34})(\la_1-\la)-f_{35}}{f_{55}}\partial_{x_5}\\
X_3(\la)&=(\la_1-\la)\partial_{x_3}+\left(1-\frac{f_{45}}{f_{55}}(\la_1-\la)\right)\partial_{x_4}+
\frac{(f_{44}-f_{35})(\la_1-\la)-f_{45}}{f_{55}}\partial_{x_5}.
\end{align*}
It is interesting that in this case the Kronecker web induced by the distribution $D_\la=\langle X_1(\la),X_2(\la),X_3(\la)\rangle$ is obviously divergence free with respect to the volume form $\om=(\la_1-\la)^2\be_\la\wedge dx_1\wedge dx_2\wedge dx_3=f_{55}dx_1\wedge dx_2\wedge dx_3\wedge dx_4\wedge dx_5$.

Alternatively, {\color{black} following \cite{bogdanov} one can achieve nontriviality of the distribution $D_\la$ by requiring vanishing of $\wedge^k\be_\la$ for some $k>2$. For instance, in the case of dimension $2n$ the condition $\wedge^n\be_\la=0$ is equivalent to vanishing
of the pfaffian $\Pf$ of the matrix of the 2-form $\be_\la$ (equivalently $\be_\infty$). In this case   rank of the 2-form $\be_\la$ is $2n-2$ and rank of $D_\la$ is 2. For example}, in dimension six, taking $N$ given by $\d_{x_i}\mapsto \la_i\d_{x_i}$ we get the equation
\begin{equation}
\label{io}\equ\begin{split}
\big((\la_5-\la_6)(\la_3-\la_4)f_{56}f_{34}-(\la_4-\la_6)(\la_3-\la_5)f_{46}f_{35}
+(\la_3-\la_6)(\la_4-\la_5)f_{36}f_{45}\big)(\la_1-\la_2)f_{12}\\
-\big((\la_5-\la_6)(\la_2-\la_4)f_{56}f_{24}-(\la_4-\la_6)(\la_2-\la_5)f_{46}f_{25}
+(\la_2-\la_6)(\la_4-\la_5)f_{26}f_{45}\big)(\la_1-\la_3)f_{13}\\
+\big((\la_5-\la_6)(\la_2-\la_3)f_{56}f_{23}-(\la_3-\la_6)(\la_2-\la_5)f_{36}f_{25}
+(\la_2-\la_6)(\la_3-\la_5)f_{26}f_{35}\big)(\la_1-\la_4)f_{14}\\
-\big((\la_4-\la_6)(\la_2-\la_3)f_{46}f_{23}-(\la_3-\la_6)(\la_2-\la_4)f_{36}f_{24}
+(\la_2-\la_6)(\la_3-\la_4)f_{26}f_{34}\big)(\la_1-\la_5)f_{15}\\
+\big((\la_2-\la_3)(\la_4-\la_5)f_{23}f_{45}-(\la_2-\la_4)(\la_3-\la_5)f_{24}f_{35}
+(\la_2-\la_5)(\la_3-\la_4)f_{25}f_{34}\big)(\la_1-\la_6)f_{16}\\
=0
\end{split}
\end{equation}
and the following parameter depending commuting vector fields spanning $D_\la$:
\begin{align*}
X_1(\la)=(\la_1-\la)\partial_{x_1}\\
+C\big[-((\la_1-\la_6)(\la_4-\la_5)f_{16}f_{45}-(\la_1-\la_5)(\la_4-\la_6)f_{15}f_{46}+
(\la_1-\la_4)(\la_5-\la_6)f_{14}f_{56}){(\la_3-\la)}\partial_{ x_3}\\
+((\la_1-\la_6)(\la_3-\la_5)f_{16}f_{35}-(\la_1-\la_5)(\la_3-\la_6)f_{15}f_{36}+
(\la_1-\la_3)(\la_5-\la_6)f_{13}f_{56}){(\la_4-\la)}\partial_{ x_4}\\
-((\la_1-\la_6)(\la_3-\la_4)f_{16}f_{34}-(\la_1-\la_4)(\la_3-\la_6)f_{14}f_{36}+
(\la_1-\la_3)(\la_4-\la_6)f_{13}f_{46}){(\la_5-\la)}\partial_{ x_5}\\
+((\la_1-\la_5)(\la_3-\la_4)f_{15}f_{34}-(\la_1-\la_4)(\la_3-\la_5)f_{14}f_{35}+
(\la_1-\la_3)(\la_4-\la_5)f_{13}f_{45}){(\la_6-\la)}\partial_{ x_6}\big]\\
X_2(\la)=(\la_2-\la)\partial_{x_2}\\
+C\big[-((\la_2-\la_4)(\la_5-\la_6)f_{24}f_{56}-(\la_4-\la_6)(\la_2-\la_5)f_{46}f_{25}+
(\la_2-\la_6)(\la_4-\la_5)f_{26}f_{45}){(\la_3-\la)}\partial_{ x_3}\\
+((\la_2-\la_3)(\la_5-\la_6)f_{23}f_{56}-(\la_3-\la_6)(\la_2-\la_5)f_{36}f_{25}+
(\la_2-\la_6)(\la_3-\la_5)f_{26}f_{35}){(\la_4-\la)}\partial_{ x_4}\\
-((\la_2-\la_3)(\la_4-\la_6)f_{23}f_{46}-(\la_3-\la_6)(\la_2-\la_4)f_{36}f_{24}+
(\la_2-\la_6)(\la_3-\la_4)f_{26}f_{34}){(\la_5-\la)}\partial_{ x_5}\\
+((\la_2-\la_3)(\la_4-\la_5)f_{23}f_{45}-(\la_3-\la_5)(\la_2-\la_4)f_{35}f_{24}+
(\la_2-\la_5)(\la_3-\la_4)f_{25}f_{34}){(\la_6-\la)}\partial_{ x_6}\big],
 \end{align*}
where $C:=\big((\la_3-\la_4)(\la_5-\la_6)f_{34}f_{56}-
(\la_3-\la_5)(\la_4-\la_6)f_{35}f_{46}+
(\la_3-\la_6)(\la_4-\la_5)f_{36}f_{45}\big)^{-1}$.

\section{Concluding remarks}
\label{concl}

\begin{enumerate}
\item {\color{black}As we mentioned in Introduction any self-dual Ricci-flat metric can be obtained by means of the Mason--Newman formalism from solutions of Plebański I equation. A natural question is whether solutions of equations $\mathbf{A_4}$, $\mathbf{B_4}$, $\mathbf{C_3}$, $\mathbf{D_3}$, and $\mathbf{E_1}$, i.e. those corresponding to the normal forms of the Nijenhuis operators with constant eigenvalues, also produce all the self-dual Ricci-flat metrics.

    This question is intemately related with the following problem of realization of Kronecker webs: given a Kronecker web of 3-web type (Definition \ref{3webdefi}) on a 4-dimensional manifold defined by means of a PNO $N:T\F_\infty\to TM$, is it true that locally there exists an extension $\widetilde{ N}:TM\to TM$ of $N$ to any of the normal forms of the Nijenhuis operators listed in Appendix? If the answer to the problem of realization is affirmative for the above mentioned cases $\mathbf{A_4}$,\ldots , $\mathbf{E_1}$, then by Corollary \ref{corrop} any pair of Mason--Newman vector fields and, consequently, any self-dual Ricci-flat metric can be obtained from solutions of the corresponding heavenly equation.

    By \cite{konopSchiefszer} any  pair of Mason--Newman vector fields can be realized in frames of $\mathbf{A_4}$-formalism (cf. Remark \ref{rem1}). The realization problem in other cases is open, however the results of \cite{pHirota}, \cite{pKronwebs}, where the problem of realization was solved for the PNOs defining Veronese webs in 3D, allow to conjecture that the answer is positive also in remaining cases.}
    {\color{black}
\item The 2-form $\be_\la$ corresponding to the general heavely equation ($\mathbf{A_4}$-case) as well as its straightforward higher-dimensional generalizations, in particular equation (\ref{io}), was first considered in \cite{bogdanov}.
\item Higher-dimensional generalizations of the I and II Plebański heavenly equations were considered in \cite{takasaki}, \cite{dunMason2000},  \cite{krynski2016}.}
\item By  \cite[Theorem 1.2]{BerjFerKrNov} all  equations of $\mathbf{A}$--$\mathbf{H}$ series are of the Monge-Amp\`{e}re type.

\item Six-dimensional equation (\ref{io}) and its  higher $2n$-dimensional generalizations $\Pf(\be_\la)=0$, where the 2-form $\be_\la$ is built by the Nijenhuis operator
    $$
    \d_{x_i}\mapsto \la_i\d_{x_i},
    $$
    generate \emph{degenerate} quadratic form \cite[(1.2)]{BerjFerKrNov}
    \begin{equation}\equ\label{qf}
    \sum_{i\le j}\frac{\d F}{\d f_{ij}}p_i p_j,
    \end{equation}
where $F(x,f,Df,D^2f)=0$ is the equation on $f$ equivalent to $\Pf(\be_\la)=0$. Indeed, it is easy to see that the matrix $\frac{\d F}{\d f_{ij}}$ is equal to
$$
\frac1{\Pf(\ad L(H))}\ad L(\widetilde{\ad L(H)}),
$$
and the equation $\Pf(\be_\la)=0$ is equivalent to $\Pf(\ad L(H))=0$,
where $H=[f_{ij}]$ is the hessian of $f$, $\widetilde{ U}$ stands for the algebraic adjoint to a matrix $U$, $L=\diag(\frac1{\la_1-\la},\ldots,\frac1{\la_{2n}-\la})$, and $\ad L(U)=LU-UL$. Since $\ad L(H)\ \widetilde{\ad L(H)}=\Pf^2(\ad L(H))I_{2n}$, we have
$$
\ad L(\widetilde{\ad L(H)})=-(\ad L(H))^{-1}\ad^2 L(H)\ \widetilde{\ad L(H)}
$$
and
$$
\det(\widetilde{\ad L(H)}))=(\Pf(\ad L(H)))^{4n-2},
$$
whence
$$
\det\left(\frac1{\Pf(\ad L(H))}\ad L(\widetilde{\ad L(H)})\right)=(\Pf(\ad L(H)))^{2n-4}\det(\ad^2 L(H)).
$$
Thus for $n>2$ with the account of the equation $\Pf(\ad L(H))=0$ quadratic form (\ref{qf}) is degenerate.

On the other hand, it is reasonable to conjecture that in dimension $2n=6$ this form has rank 4 on the solutions of $\Pf(\ad L(H))=0$, which by  \cite[Theorem 1.3]{BerjFerKrNov} would imply the Monge-Amp\`{e}re property for (\ref{io}). {\color{black} And indeed, one checks directly that the LHS of (\ref{io}) is
$$
-\sum_{i>j>k=1}^6p_{ijk}N_{ijk}\la_i\la_j\la_k,
$$
where $N_{ijk}$ is a determinant of the submatrix of the Hessian $\mathrm{Hess}(f)$ obtained by picking the columns with the numbers $i>j>k$ and the lines with the numbers  belonging to the naturally ordered set $c(i,j,k):=\{1,\ldots,6\}\setminus\{i,j,k\}$ and $p_{ijk}$ is the parity of the permutation ${\color{black}(1,\ldots,6)\mapsto}(i,j,k,c(i,j,k))$.}
\item Equation $\mathbf{A_0}$ appeared in \cite{BerjFerKrNov} and it was shown that $\mathbf{A_0}$ and $\mathbf{A_4}$ (which is the same as (\ref{schief})) are contact nonequivalent.
\item Equation $\mathbf{G_0}$ fits to the definition of equations of ``first heavenly type'' \cite{BerjFerKrNov} and by results of that paper should be contact equivalent to one of three equations listed in \cite[Table 2]{BerjFerKrNov}.
\item Equation $\mathbf{D_3}$ appeared in \cite[(168)]{konopSchiefszer} as a reduction of the so-called TED equation.
\item As we noticed in Section \ref{sec1}, with any classical 3-web one can associate the canonical Chern connection \cite{chern}. On the othe hand, given a Kronecker web $\{\F_\la\}$ of 3-web type, one can construct a classical 3-web $\F_{\la_1},\F_{\la_2},\F_{\la_3}$ by picking up any three values $\la_i$ of the parameter $\la$. The corresponding Chern connection will be torsionless \cite{Nagy} and, moreover, in the divergence free case in 4D one can choose $\la_i$  in such a way that this connection coincide up to a constant factor with the Levi-Civita connection of the metric (\ref{intro2}). This allows to derive different properties of the Levi-Civita connection from those of the Chern connection.
\item An intriguing task is to characterize in terms of the Chern connection those 4-dimensional divergence free Kronecker webs of 3-web type which can be included in a Veronese web (see our discussion of relations between heavenly PDEs and Hirota dispesionless systems in Introduction).
\end{enumerate}

\section{Appendix: Normal forms of Nijenhuis operators in 4D }\label{app}

In papers \cite{t4,grifoneMehdi} the authors obtained a local classification of Nijenhuis operators $N:TM\to TM$ (in a vicinity of a regular point \cite[p. 451]{t4})
under the additional assumption of existence of a complete family of the conservation laws. This assumption is equivalent to vanishing of
the invariant $P_N$, which is automatically trivial in the case of cyclic $N$ \cite[p. 450]{t4}, i.e., when the space $T_xM$ is cyclic for $N_x$ for any $x\in M$.
Here we recall these normal forms obtained in this case for 4-dimensional $M$ {\color{black} in complex category} (cf. also \cite{bolsKonMatBIG}) (cases $\mathbf{A}$--$\mathbf{E}$). We also list some obvious pairwise nonequivalent noncyclic normal forms of Nijenhuis operators (cases $\mathbf{F}$--$\mathbf{H}$).

Below  $(x_1,x_2,x_3,x_4)$ denotes a system of local coordinates
 and  $\la_1,\la_2,\la_3,\la_4$ are pairwise distinct constants.
In the the follwing list   $N$ stays for the matrix of a Nijenhuis operator in the basis $\{\partial_{x_i}\}$ and we also indicate the corresponding
Jordan normal form $J$.
\begin{itemize}
\item[\bf ${\bold A_0}$.] $N_{A0}=N_{A0}(x_1,x_2,x_3,x_4):=\left[
           \begin{array}{cccc}
             x_1  & 0 & 0 & 0 \\
             0 & x_2 & 0 & 0\\
             0 & 0 & x_3 & 0\\
              0 & 0 & 0 & x_4
           \end{array}
         \right]$, \quad $J_{A0}=N_{A0}$.\\
\item[\bf ${\bold A_1}$.] $N_{A1}:=N_{A0}(x_1,x_2,x_3,\la_4)$.
\item[\bf ${\bold A_2}$.] $N_{A2}:=N_{A0}(x_1,x_2,\la_3,\la_4)$.
\item[\bf ${\bold A_3}$.] $N_{A3}:=N_{A0}(x_1,\la_2,\la_3,\la_4)$.
\item[\bf ${\bold A_4}$.] $N_{A4}:=N_{A0}(\la_1,\la_2,\la_3,\la_4)$.
\item[\bf ${\bold B_0}$.] $N_{B0}=N_{B0}(x_1,x_3,x_4):=\left[
           \begin{array}{cccc}
             x_1  & 0 & 0 & 0 \\
             1 & x_1 & 0 & 0\\
             0 & 0 & x_3 & 0\\
              0 & 0 & 0 & x_4
           \end{array}
         \right]$, \quad  $J_{B0}=N_{B0}$.
\item[\bf ${\bold B_1}$.] $N_{B1}:=N_{B0}(x_1,x_3,\la_4)$.
\item[\bf ${\bold B_2}$.] $N_{B2}:=N_{B0}(x_1,\la_3,\la_4)$.
\item[\bf ${\bold B_3}$.] $N_{B3}:=N_{B0}(\la_1,x_3,\la_4)$.
\item[\bf ${\bold B_4}$.] $N_{B4}:=N_{B0}(\la_1,\la_3,\la_4)$.
\item[\bf ${\bold C_0}$.] $N_{C0}=N_{C0}(x_2,x_3,x_4):=\left[
           \begin{array}{cccc}
             x_3  & 0 & 1 & 0 \\
             1 & x_3 & -x_2 & 0\\
             0 & 0 & x_3 & 0\\
              0 & 0 & 0 & x_4
           \end{array}
         \right]$, \quad
         $J_{C0}=J_{C0}(x_3,x_4):=\left[
           \begin{array}{cccc}
             x_3  & 0 & 0 & 0 \\
             1 & x_3 & 0 & 0\\
             0 & 1 & x_3 & 0\\
              0 & 0 & 0 & x_4
           \end{array}
         \right]$.
\item[\bf ${\bold C_1}$.] $N_{C1}:=N_{C0}(x_2,x_3,\la_4)$.
\item[\bf ${\bold C_2}$.] $N_{C2}:=J_{C0}(\la_3,x_4)$.
\item[\bf ${\bold C_3}$.] $N_{C3}:=J_{C0}(\la_3,\la_4)$.

\item[\bf ${\bold D_0}$.] $N_{D0}=N_{D0}(x_1,x_3):=\left[
           \begin{array}{cccc}
             x_1  & 0 & 0 & 0 \\
             1 & x_1 & 0 & 0\\
             0 & 0 & x_3 & 0\\
              0 & 0 & 1 & x_3
           \end{array}
         \right]$, \quad
         $J_{D0}:=N_{D0}$.
\item[\bf ${\bold D_1}$.] $N_{D1}:=N_{D0}(x_1,\la_3)$.
\item[\bf ${\bold D_2}$.] $N_{D2}:=N_{D0}(\la_1,x_3)$.
\item[\bf ${\bold D_3}$.] $N_{D3}:=N_{D0}(\la_1,\la_3)$.
\item[\bf ${\bold E_0}$.] $N_{E0}=N_{E0}(x_2,x_3,x_4):=\left[
           \begin{array}{cccc}
             x_4  & 0 & 0 & 1 \\
             1 & x_4 & 0 & -x_2\\
             0 & 1 & x_4 & -2x_3\\
              0 & 0 & 0 & x_4
           \end{array}
         \right]$, \quad
         $J_{E0}=J_{E0}(x_4):=\left[
           \begin{array}{cccc}
             x_4  & 0 & 0 & 0 \\
             1 & x_4 & 0 & 0\\
             0 & 1 & x_4 & 0\\
              0 & 0 & 1 & x_4
           \end{array}
         \right]$.
\item[\bf ${\bold E_1}$.] $N_{E1}:=J_{E0}(\la_4)$.
\item[\bf ${\bold F_0}$.] $N_{F0}=N_{F0}(x_1,x_3,x_4):=\left[
           \begin{array}{cccc}
             x_1  & 0 & 0 & 0 \\
             0 & x_1 & 0 & 0\\
             0 & 0 & x_3 & 0\\
              0 & 0 & 0 & x_4
           \end{array}
         \right]$, \quad $J_{F0}=N_{F0}$.\\
\item[\bf ${\bold F_1}$.] $N_{F1}:=N_{F0}(x_1,x_3,\la_4)$.
\item[\bf ${\bold F_2}$.] $N_{F2}:=N_{F0}(x_1,\la_3,\la_4)$.
\item[\bf ${\bold F_3}$.] $N_{F3}:=N_{F0}(\la_1,x_3,x_4)$.
\item[\bf ${\bold F_4}$.] $N_{F4}:=N_{F0}(\la_1,\la_3,x_4)$.
\item[\bf ${\bold F_5}$.] $N_{F5}:=N_{F0}(\la_1,\la_3,\la_4)$.
\item[\bf ${\bold G_0}$.] $N_{G0}=N_{G0}(x_1,x_3):=\left[
           \begin{array}{cccc}
             x_1  & 0 & 0 & 0 \\
             0 & x_1 & 0 & 0\\
             0 & 0 & x_3 & 0\\
              0 & 0 & 0 & x_3
           \end{array}
         \right]$, \quad $J_{G0}=N_{G0}$.\\
\item[\bf ${\bold G_1}$.] $N_{G1}:=N_{G0}(x_1,\la_3)$.
\item[\bf ${\bold G_2}$.] $N_{G2}:=N_{G0}(\la_1,x_3)$.
\item[\bf ${\bold G_3}$.] $N_{G3}:=N_{G0}(\la_1,\la_3)$.
\item[\bf ${\bold H_0}$.] $N_{H0}=N_{H0}(x_1,x_3):=\left[
           \begin{array}{cccc}
             x_1  & 0 & 0 & 0 \\
             1 & x_1 & 0 & 0\\
             0 & 0 & x_3 & 0\\
              0 & 0 & 0 & x_3
           \end{array}
         \right]$, \quad $J_{H0}=N_{H0}$.\\
\item[\bf ${\bold H_1}$.] $N_{H1}:=N_{H0}(x_1,\la_3)$.
\item[\bf ${\bold H_2}$.] $N_{H2}:=N_{H0}(\la_1,x_3)$.
\item[\bf ${\bold H_3}$.] $N_{H3}:=N_{H0}(\la_1,\la_3)$.
\end{itemize}

{\color{black}
\section{Acknowledgements}

The authors are indebted to Boris Kruglikov for useful comments {\color{black} and to Wolfgang Schief for indicating  reference \cite{bogdanov}.}}


\begin{thebibliography}{BFKN20}

\bibitem[AS02]{martinezAlShabat}
L.~Mart\'{i}nez Alonso and A.B. Shabat, \emph{Energy-dependent potentials
  revisited: a universal hierarchy of hydrodynamic type}, Phys. Lett. A
  \textbf{300} (2002), 58--64.

\bibitem[BFKN20]{BerjFerKrNov}
S.~Berjawi, E.~V. Ferapontov, B.~Kruglikov, and V.~Novikov, \emph{Second-order
  {PDE}s in four dimensions with half-flat conformal structure}, Proc. R. Soc.
  A \textbf{476} (2020), 20190642.

\bibitem[Bog15]{bogdanov}
L. ~V. Bogdanov,
\emph{Doubrov–Ferapontov general heavenly equation and the hyper-Kähler hierarchy},
 J. Phys. A: Math. Theor. \textbf{48} (2015), 235202.


\bibitem[BKM21]{bolsKonMat}
A.~Bolsinov, A.~Konyaev, and V.~Matveev, \emph{Applications of {N}ijenhuis
  geometry: non-degenerate singular points of {P}oisson--{N}ijenhuis
  structures}, Europ. J. Math. \textbf{8} (2021), 1355–1376.

\bibitem[BKM22]{bolsKonMatBIG}
\bysame, \emph{{N}ijenhuis geometry}, Adv. Math. \textbf{394} (2022), 108001.

\bibitem[Che36]{chern}
S.~S. Chern, \emph{Einen {I}nvariantentheorie der 3-gewebe aus r-dimensionalen
  {M}annifaltigkeiten in $\mathbb{R}_{2r}$}, Abh. Math. Sem. Univ. Hamburg
  \textbf{11} (1936), 333--358.

\bibitem[DZ23]{domZub}
W.~Domitrz and M.~Zubilewicz, \emph{Local invariants of divergence-free webs},
  Anal. Math. Phys. \textbf{13} (2023).


\bibitem[DF10]{doubrovFerapontov}
B.~Doubrov and E.~V. Ferapontov, \emph{On the integrability of symplectic
  {M}onge-{A}mp\'{e}re equations}, J. Geom. Phys. \textbf{60} (2010), no.~10,
  1604--1616.

\bibitem[Dun04]{dunajski}
M.~Dunajski, \emph{A class of {E}instein-{W}eyl spaces associated to an
  integrable system of hydrody- namic type}, J. Geom. Phys. \textbf{51} (2004),
  126--137.

\bibitem[Dun10]{dunajskiBook}
\bysame, \emph{Solitons, instantons, and twistors}, Oxford {U}niversity
  {P}ress, 2010.


\bibitem[DK14]{dunajskiKr}
M.~Dunajski and W.~Kry\'{n}ski, \emph{Einstein-{W}eyl geometry, dispersionless
  {H}irota equation and {V}eronese webs}, Math. Proc. Camb. Phil. Soc.
  \textbf{157} (2014), 139--150.

\bibitem[DM00]{dunMason2000}
M.~Dunajski and L. ~J. ~Mason, \emph{Hyper-{K}\"{a}hler hierarchies and their twistor theory}, Comm. Math. Phys. \textbf{213} (2000), 641--672.



\bibitem[FK21]{ferapontovKruglikov}
E.~Ferapontov and B.~Kruglikov, \emph{Dispersionless integrable hierarchies and
  {GL(2, R)} geometry}, Math. Proc. Camb. Phil. Soc \textbf{170} (2021),
  129--154.

\bibitem[FP79]{FinPleb}
J.~D. Finley and J.~F. Pleba\'nski, \emph{The classification of all
  $\mathcal{H}$ spaces admitting a {K}illing vector}, J. Math. Phys.
  \textbf{20} (1979), 1938--1945.

\bibitem[Gin82]{gindikin82}
S.~Gindikin, \emph{Bundles of differential forms and the {E}instein equation},
  Nuclear Phys. (in Russian) \textbf{36} (1982), no.~2(8), 537--548.

\bibitem[Gin85]{gindikin85}
S.~Gindikin, \emph{Some solutions of the selfdual Einstein equations}, Funktsional. Anal. i Prilozhen. (in Russian) \textbf{19}  (1985), no. ~3, 58--60.

\bibitem[GM97]{grifoneMehdi}
J.~Grifone and M.~Mehdi, \emph{Existence of conservation laws and
  characterization of recursion operators for completely integrable systems},
  Trans. Amer. Math. Soc. \textbf{349} (1997), 4609--4633.

\bibitem[Gra93]{grant}
J.~Grant, \emph{On self-dual gravity}, Phys. Rev. D \textbf{48} (1993),
  2606--2612.

\bibitem[GZ91]{gz2}
I.~Gelfand and I.~Zakharevich, \emph{Webs, {V}eronese curves, and bihamiltonian
  systems}, J. Funct. Anal. \textbf{99} (1991), 150--178.

\bibitem[Hus94]{husain}
V.~Husain, \emph{Self-dual gravity and the chiral model}, Phys. Rev. Lett.
  \textbf{72} (1994), 800.

  \bibitem[KSS21]{konopSchiefszer}
B.~Konopelchenko, W.~Schief, and A.~Szereszewski, \emph{Self-dual {E}instein
  spaces and the general heavenly equation. {E}igenfunctions as coordinates},
  Class. Quantum Grav. \textbf{38} (2021), 045007.

\bibitem[Kos85]{koszul}
J.-L. Koszul, \emph{Crochets de {S}chouten {N}ijenhuis et cohomology},
  Ast\'{e}risque, Sc. Math. de france, hors s\'{e}rie, 1985, 257--271.





\bibitem[KSM90]{mks}
Y.~Kosmann-Schwarzbach and F.~Magri, \emph{Poisson--{N}ijenhuis structures},
  Ann. Inst. Henri Poincar\'{e} \textbf{53} (1990), 35--81.

\bibitem[KP17]{pHirota}
B.~Kruglikov and A.~Panasyuk, \emph{Veronese webs and nonlinear {PDE}s}, J.
  Geom. Phys. \textbf{115} (2017), 45--60.

\bibitem[Kry16]{krynski2016} W. ~Kryński, \emph{Webs and the {P}lebański equation}, Math. Proc. Camb. Phil.
Soc. \textbf{161(3) }(2016) 455--468.

\bibitem[Kry18]{KrynskiDef}
W.~Kry\'{n}ski, \emph{On deformations of the dispersionless {H}irota equation},
  J. Geom. Phys. \textbf{127} (2018), 46--54.

\bibitem[MN89]{masonNewman}
L.~J. Mason and E.T. Newman, \emph{A connection between the {E}insteinand
  {Y}ang--{M}ills equations}, Commun. Math. Phys. \textbf{121} (1989),
  659--668.

\bibitem[MW96]{masonWoodhouse}
L.~J. Mason and N.~M.~J. Woodhouse, \emph{Integrability, {S}elf-{D}uality and
  {T}wistor {T}heory}, LMS Monograph New Series, vol.~15, Oxford, 1996.

\bibitem[Nag01]{Nagy}
P.~T. Nagy, \emph{Webs and curvature}, Web Theory and Related Topics, World
  Scientific, 2001, pp.~48--91.

\bibitem[Pan19]{pKronwebs}
A.~Panasyuk, \emph{Kronecker webs, {N}ijenhuis operators, and nonlinear
  {PDE}s}, Banach Center Publications \textbf{117} (2019), 177--210.

\bibitem[Par92]{park}
Q.~H. Park, \emph{2d sigma model approach to 4d instantons}, Int. J. Mod. Phys.
  \textbf{A7} (1992), 1415.

\bibitem[Pen76]{penrose}
R.~Penrose, \emph{Nonlinear gravitons and curved twistor theory}, Gen. Rel.
  Grav. \textbf{7} (1976), 31--52.

\bibitem[Ple75]{plebanski}
J.~F. Pleba\'{n}ski, \emph{Some solutions of complex {E}instein equations}, J.
  Math. Phys. \textbf{16} (1975), 2395--2402.

\bibitem[Rig96]{rigal}
M.-H. Rigal, \emph{Geometrie globale des systemes bihamiltoniens en dimension
  impaire}, Ph.D. thesis, l'Universit\'{e} Montpellier II, 1996.

\bibitem[Sch96]{schief}
W.~K. Schief, \emph{Self-dual {E}instein spaces via a permutability theorem for
  the {T}zitzeica equation}, Phys. Lett. A \textbf{223} (1996), 55--62.

\bibitem[Tab93]{taba}
S.~Tabachnikov, \emph{Geometry of {L}agrangian and {L}egendrian 2-web}, Diff.
  Geometry Appl. \textbf{3} (1993), 265--284.

\bibitem[Tak89]{takasaki}
K. ~Takasaki, \emph{An infinite number of hidden variables in
hyper-Kähler metrics}, J. Math. Phys. \textbf{30} (1989), 1515--1521.

\bibitem[Tur96]{t4}
F.-J. Turiel, \emph{Classification of (1,1) tensor fields and bihamiltonian
  structures}, Banach Center Publications \textbf{33} (1996), 449--458.

\bibitem[Tur97]{turielDDN}
\bysame, \emph{L'\'{e}quation $2i\partial\bar{\partial}f=\beta$ pour les
  tenseurs de {N}ijenhuis de type $(1,1)$}, C. R. Acad. Sci. Paris, S\'{e}rie I
  \textbf{325} (1997), 1313--1316.

\bibitem[Zak01]{z1}
I.~Zakharevich, \emph{Kronecker webs, bihamiltonian structures, and the method
  of argument translation}, Transform. Groups \textbf{6} (2001), 267--300.

\end{thebibliography}
\providecommand{\bysame}{\leavevmode\hbox to3em{\hrulefill}\thinspace}
\providecommand{\MR}{\relax\ifhmode\unskip\space\fi MR }
\providecommand{\MRhref}[2]{%
  \href{http://www.ams.org/mathscinet-getitem?mr=#1}{#2}
}
\providecommand{\href}[2]{#2}

\end{document}